\documentclass[a4paper,10pt,twoside]{article}

\usepackage[utf8]{inputenc}
\usepackage{amsmath, amssymb, amsthm}
\usepackage{amstext, amsfonts, a4}
\usepackage{dsfont}
\usepackage{latexsym}
\usepackage{mathrsfs}
\usepackage[all,cmtip]{xy}
\usepackage{graphicx}
\usepackage{enumerate}
\usepackage{hyperref,color}
\usepackage{stmaryrd}
\usepackage[shortlabels]{enumitem}

\usepackage{url}
\usepackage[pdftex,a4paper,left=3cm,right=3cm,top=3cm,bottom=3cm]{geometry}

\def\ra{\rightarrow}

\def\lra{\longrightarrow}

\def\lmapsto{\longmapsto}

\def\sA{\mathscr{A}}

 \def\sZ{\mathscr{Z}}

\def\bbC{\mathbb{C}}
\def\bbF{\mathbb{F}}\def\bbG{\mathbb{G}}

\def\bbP{\mathbb{P}}
\def\bbQ{\mathbb{Q}}\def\bbT{\mathbb{T}}

\def\bbZ{\mathbb{Z}}

\def\cB{\mathcal{B}}
\def\cF{\mathcal{F}}\def\cH{\mathcal{H}}
\def\cI{\mathcal{I}}
\def\cM{\mathcal{M}}\def\cO{\mathcal{O}}\def\cP{\mathcal{P}}

\def\cW{\mathcal{W}}
\def\cZ{\mathcal{Z}}

\def\bfB{\mathbf{B}}
\def\bfG{\mathbf{G}}

\def\bfT{\mathbf{T}}

\def\bfq{\mathbf{q}}

\def\fS{\mathfrak{S}}

\DeclareMathOperator{\aff}{aff}
\DeclareMathOperator{\Aut}{Aut}

\DeclareMathOperator{\diag}{diag}

\DeclareMathOperator{\End}{End}

\DeclareMathOperator{\Hom}{Hom}

\DeclareMathOperator{\id}{id}
\DeclareMathOperator{\Id}{Id}

\DeclareMathOperator{\ind}{ind}

\DeclareMathOperator{\nil}{nil}

\DeclareMathOperator{\perm}{perm}

\DeclareMathOperator{\Pic}{Pic}

\DeclareMathOperator{\pt}{pt}
\DeclareMathOperator{\proj}{proj}

\DeclareMathOperator{\Spec}{Spec}

\DeclareMathOperator{\Sym}{Sym}

\newtheorem{counter}[subsection]{$\!\!$}

\newenvironment{Prop}{\begin{counter} {\bf Proposition.}}{\end{counter}}
\newenvironment{Lem}{\begin{counter} {\bf Lemma.}}{\end{counter}}

\newenvironment{Pt}{\begin{counter} \rm}{\end{counter}}

\newenvironment{Proof}{{\flushleft \bf Proof :}}{\hfill $\square$ \vspace{5mm}}

\newtheorem{counter*}[subsubsection]{$\!\!$}
\newenvironment{Def*}{\begin{counter*} {\bf Definition.}}{\end{counter*}}
\newenvironment{Not*}{\begin{counter*} \rm {\bf Notation.}}{\end{counter*}}
\newenvironment{Notss*}{\begin{counter*} \rm {\bf Notations.}}{\end{counter*}}
\newenvironment{DefNot*}{\begin{counter*} \rm {\bf Definition-Notation.}}{\end{counter*}}
\newenvironment{Nots*}{\begin{counter*} \rm {\bf Notations.}}{\end{counter*}}
\newenvironment{Prop*}{\begin{counter*} {\bf Proposition.}}{\end{counter*}}
\newenvironment{Lem*}{\begin{counter*} {\bf Lemma.}}{\end{counter*}}
\newenvironment{Cor*}{\begin{counter*} {\bf Corollary.}}{\end{counter*}}
\newenvironment{Th*}{\begin{counter*} {\bf Theorem.}}{\end{counter*}}
\newenvironment{Rem*}{\begin{counter*} \rm {\bf Remark.}}{\end{counter*}}
\newenvironment{Ex*}{\begin{counter*} \rm {\bf Example.}}{\end{counter*}}
\newenvironment{Exs*}{\begin{counter*} \rm {\bf Examples.}}{\end{counter*}}
\newenvironment{Pt*}{\begin{counter*} \rm}{\end{counter*}}
\newenvironment{Q*}{\begin{counter*} \rm {\bf Question.}}{\end{counter*}}

\title{\textbf{\huge{Mod $p$ Hecke algebras \\ and dual equivariant cohomology I: \\ the case of $GL_2$}}}

\author{Cédric PEPIN and Tobias SCHMIDT}
\date{\today}

\begin{document}

\maketitle

\begin{abstract}
Let $F$ be a $p$-adic local field and $\bfG=GL_2$ over $F$. Let $\cH^{(1)}$ be the pro-$p$ Iwahori-Hecke algebra of the group $\bfG(F)$ with coefficients in the algebraic closure $\overline{\bbF}_p$. We show that the supersingular irreducible $\cH^{(1)}$-modules of dimension $2$ can be realized through
the equivariant cohomology of the flag variety of the Langlands dual group $\widehat{\bfG}$ over $\overline{\bbF}_p$.
\end{abstract}

\tableofcontents

\section{Introduction}
Let $F$ be a finite extension of $\bbQ_p$ with residue field $\bbF_q$ and let $\bfG$ be a connected split reductive group over $F$. Let $\cH=R[I\setminus\bfG(F)/I]$ be the Iwahori-Hecke algebra associated to an Iwahori subgroup $I\subset \bfG(F)$, with coefficients in an algebraically closed field $R$. On the other hand, let $\widehat{\bfG}$ be the Langlands dual group of $\bfG$ over $R$, and $\widehat{\cB}$ the flag variety of Borel subgroups of $\widehat{\bfG}$ over $R$.

\vskip5pt

When $R=\bbC$, the irreducible $\cH$-modules appear as subquotients of the Grothendieck group $K^{\widehat{\bfG}}(\widehat{\cB})_{\bbC}$ of
$\widehat{\bfG}$-equivariant coherent sheaves on $\widehat{\cB}$. As such they can be parametrized by the isomorphism classes of irreducible tame $\widehat{\bfG}(\bbC)$-representations of the absolute Galois
group ${\rm Gal}(\overline{F}/F)$ of $F$, thereby realizing the tame local Langlands correspondence (in this setting also called the Deligne-Lusztig conjecture for Hecke modules): Kazhdan-Lusztig \cite{KL87}, Ginzburg \cite{CG97}. The idea of studying various cohomological invariants of the flag variety by means of Hecke operators (nowadays called Demazure operators) goes back to earlier work of Demazure \cite{D73,D74}. The approach to the Deligne-Lusztig conjecture is based on the construction of a natural $\cH$-action on the whole $K$-group $K^{\widehat{\bfG}}(\widehat{\cB})_{\bbC}$ which identifies the center of $\cH$ with the $K$-group of the base point $K^{\widehat{\bfG}}({\rm pt})_{\bbC}$. The finite part of $\cH$ acts thereby via appropriate $\bfq$-deformations of Demazure operators.

\vskip5pt

When $R=\overline{\bbF}_q$ any irreducible
$\widehat{\bfG}(\overline{\bbF}_q)$-representation of
${\rm Gal}(\overline{F}/F)$ is tame and the Iwahori-Hecke algebra needs to be replaced by the bigger pro-$p$-Iwahori-Hecke algebra
$$\cH^{(1)}=\overline{\bbF}_q[I^{(1)}\setminus \bfG(F)/I^{(1)}].$$ Here, $I^{(1)}\subset I$ is the unique pro-$p$ Sylow subgroup of $I$.
The algebra $\cH^{(1)}$ was introduced by Vign\'eras and its structure theory developed in a series of papers
\cite{V04,V05,V06,V14,V15,V16,V17}. The class of so-called {\it supersingular} irreducible $\cH^{(1)}$-modules figures prominently among all irreducible $\cH^{(1)}$-modules, since it is expected to be related to the arithmetic over the field $F$. For $\bfG=GL_n$, there is a distinguished correspondence between supersingular irreducible $\cH^{(1)}$-modules of dimension $n$ and irreducible $GL_n(\overline{\bbF}_q)$-representations of ${\rm Gal}(\overline{F}/F)$: Breuil \cite{Br03}, Vign\'eras \cite{V04}, \cite{V05}, Colmez \cite{C10}, Grosse-Klönne \cite{GK16}, \cite{GK18}.

\vskip5pt
Our aim is to show that the supersingular irreducible $\cH^{(1)}$- modules of dimension $n$ can again be realized as subquotients of some
$\widehat{\bfG}$-equivariant cohomology theory of the flag variety $\widehat{\cB}$ over $\overline{\bbF}_q$, although in a way different from the $\bbC$-coefficient case. Here we discuss the case $n=2$, and we will treat the case of general $n$ in a subsequent article \cite{PS2}.

\vskip5pt

From now on, let $R=\overline{\bbF}_q$ and $\bfG=GL_2$.
The algebra $\cH^{(1)}$ splits as a direct product of subalgebras $\cH^{\gamma}$ indexed by the orbits $\gamma$ of $\fS_2$ in the set of characters of $(\bbF_q^\times)^2$, namely the Iwahori components corresponding to trivial orbits, and the regular components. Accordingly, the category of $\cH^{(1)}$-modules decomposes as the product of the module categories for the component algebras. In each component sits a unique supersingular
module of dimension $2$ with given central character. On the dual side, we have the projective line $\widehat{\cB}=\bbP^1_{\overline{\bbF}_q}$ over $\overline{\bbF}_q$ with its natural action by fractional transformations of the algebraic group $\widehat{\bfG}=GL_2(\overline{\bbF}_q)$.

\vskip5pt

For a non-regular orbit $\gamma$, the component algebra $\cH^{\gamma}$
is isomorphic to the mod $p$ Iwahori-Hecke algebra $\cH=\overline{\bbF}_q[I\setminus\bfG(F)/I]$ and the quadratic relations in $\cH$ are idempotent of type $T^2_s=-T_s$. The $\widehat{\bfG}$-equivariant $K$-theory $K^{\widehat{\bfG}}(\widehat{\cB})_{\overline{\bbF}_q}$ of $\widehat{\cB}$ comes with an action of the classical Demazure operator at $\bfq=0$. Our first result is that this action extends uniquely to an action of the full algebra $\cH$ on $K^{\widehat{\bfG}}(\widehat{\cB})_{\overline{\bbF}_p}$, which is faithful and which identifies the center $Z(\cH)$ of $\cH$ with the base ring
$K^{\widehat{\bfG}}({\rm pt})_{\overline{\bbF}_q}$. It is constructed from natural presentations of the algebras $\cH$ and $Z(\cH)$ \cite{V04} and through the characteristic homomorphism
$$\bbZ[\Lambda]\stackrel{\simeq}{\longrightarrow}K^{\widehat{\bfG}}(\widehat{\cB})$$ which identfies the equivariant $K$-ring with the group ring of characters $\Lambda$ of a maximal torus in $\widehat{\bfG}$. In particular, everything is explicit. We finally show that, given a supersingular central character
$\theta: Z(\cH)\ra\overline{\bbF}_q$, the central reduction
$K^{\widehat{\bfG}}(\widehat{\cB})_{\theta}$ is isomorphic to the unique supersingular $\cH$-module of dimension $2$ with central character 
$\theta$.

\vskip5pt

For a regular orbit $\gamma$, the component algebra $\cH^{\gamma}$
is isomorphic to Vign\'eras second Iwahori-Hecke algebra $\cH_2$ \cite{V04}. It can be viewed as a certain twisted version of two copies of the mod $p$ nil Hecke ring  $\cH^{\rm nil}$ (introduced over the complex numbers by Kostant-Kumar \cite{KK86}). In particular,
the quadratic relations are nilpotent of type $T_s^2=0$. The $\widehat{\bfG}$-equivariant intersection theory $CH^{\widehat{\bfG}}(\widehat{\cB})_{\overline{\bbF}_q}$ of $\widehat{\cB}$ comes with an action of the classical Demazure operator at $\bfq=0$. We show that this action extends to a faithful action of $\cH^{\rm nil}$ on $CH^{\widehat{\bfG}}(\widehat{\cB})_{\overline{\bbF}_q}$.
To incorporate the twisting, we then pass to the square $\widehat{\cB}^2$ of $\widehat{\cB}$ and extend the action to a faithful action of $\cH_2$ on  $CH^{\widehat{\bfG}}(\widehat{\cB}^2)_{\overline{\bbF}_q}$. The action identifies a large part $Z^{\circ}(\cH_2)$ of the center $Z(\cH_2)$ with the base ring $CH^{\widehat{\bfG}}({\rm pt})_{\overline{\bbF}_q}$. As a technical point, one actually has to pass to a certain localization of the Chow groups to realize these actions, but we do not go into this in the introduction. As in the non-regular case, the action is constructed from natural presentations of the algebras $\cH_2$ and $Z(\cH_2)$ \cite{V04} and through the characteristic homomorphism
$$\Sym (\Lambda) \stackrel{\simeq}{\longrightarrow} CH^{\widehat{\bfG}}(\widehat{\cB})$$ which identifies the equivariant Chow ring with the symmetric algebra on the character group $\Lambda$. So again, everything is explicit.
We finally show that, given a supersingular central character
$\theta: Z(\cH_2)\ra\overline{\bbF}_q$, the semisimplification of the $Z^{\circ}(\cH_2)$-reduction of (the localization of) $CH^{\widehat{\bfG}}(\widehat{\cB}^2)_{\overline{\bbF}_q}$ equals a direct sum of four copies of the unique supersingular $\cH_2$-module of dimension $2$ with central character $\theta$.

 \vskip5pt

 In a final section we discuss the aforementioned bijection between supersingular irreducible $\cH^{(1)}_{\overline{\bbF}_q}$-modules of dimension $2$ and irreducible smooth $GL_2(\overline{\bbF}_q)$-representations of ${\rm Gal}(\overline{F}/F)$ in the light of our geometric language.

 \vskip5pt

 {\it Notation:} In general, the letter $F$ denotes a locally compact complete non-archimedean field with ring of integers $o_F$.
 Let $\bbF_q$ be its residue field, of characteristic $p$ and cardinality $q$. We denote by $\bfG$ the algebraic group ${\rm GL}_{2}$ over $F$ and by $G:=\bfG(F)$ its group of $F$-rational points. Let $\bfT\subset\bfG$ be the torus of diagonal matrices. Finally, $I\subset G$ denotes the upper triangular standard Iwahori subgroup and $I^{(1)}\subset I$ denotes the unique pro-$p$ Sylow subgroup of $I$. Without further mentioning, all modules will be left modules.

\section{The pro-$p$-Iwahori-Hecke algebra}

Let $R$ be any commutative ring. The \emph{pro-$p$ Iwahori Hecke algebra of $G$ with coefficients in $R$} is defined to be the convolution algebra
$$
\cH_{R}^{(1)}(q):=(R[I^{(1)}\backslash G/I^{(1)}],\star)
$$
generated by the $I^{(1)}$-double cosets in $G$. In the sequel, \emph{we will assume that $R$ is an algebra over the ring}
$$
\bbZ[\frac{1}{q-1},\mu_{q-1}].
$$
The first examples we have in mind are $R=\bbF_q$ or its algebraic closure $R=\overline{\bbF}_q$.
\subsection{Weyl groups and cocharacters}
\begin{Pt*}\label{cochar}
We denote by $$\Lambda=\Hom(\bbG_m,\bfT)=\bbZ\eta_1\oplus\bbZ\eta_2\simeq\bbZ\oplus\bbZ$$
the lattice of cocharacters of $\bfT$ with standard basis $\eta_1(x)=\diag(x,1)$ and  $\eta_2(x)=\diag(1,x)$. Then $\alpha=(1,-1)\in\Lambda$ is a root and
\begin{eqnarray*}
s=s_{\alpha}=s_{(1,-1)}:\bbZ\oplus\bbZ&\lra&\bbZ\oplus\bbZ \\
(n_1,n_2)&\lmapsto& (n_2,n_1)
\end{eqnarray*}
is the associated reflection generating the Weyl group $W_0=\{1,s\}$. The element $s$ acts on $\Lambda$ and hence also on the group ring
$\bbZ[\Lambda]$. The two invariant elements
$$\xi_1:=e^{(1,0)}+e^{(0,1)}\hskip15pt \text{and} \hskip15pt \xi_2:=e^{(1,1)}$$
in $\bbZ[\Lambda]^{s}$
define a ring isomorphism
\begin{eqnarray*}
\xi^+: \bbZ[\Lambda^{+}]=\bbZ\big{[} e^{(1,0)},(e^{(1,1)})^{\pm 1}\big{]} & \stackrel{\simeq} {\longrightarrow}  & \bbZ[\Lambda]^{s} \\
e^{(1,0)} & \lmapsto & \xi_1\\
e^{(1,1)} & \lmapsto & \xi_2
\end{eqnarray*}
where $\Lambda^{+}:=\bbZ_{\geq 0}(1,0)\oplus\bbZ(1,1)$ is the monoid of dominant cocharacters.
\end{Pt*}

\begin{Pt*}

We introduce the affine Weyl group $W_{\aff}$ and the Iwahori-Weyl group $W$ of $G$:

$$
 W_{\aff}:=e^{\bbZ(1,-1)}\rtimes W_0 \subset W:=e^\Lambda\rtimes W_0.
$$
With $$u:=e^{(1,0)}s=se^{(0,1)}
$$
one has
$W= W_{\aff}\rtimes \Omega$
where $\Omega=u^{\bbZ}\simeq\bbZ.$
Let $s_0=e^{(1,-1)}s=se^{(-1,1)}=usu^{-1}$.
Recall that the pair $( W_{\aff},\{s_0,s\})$ is a Coxeter group and its length
function $\ell$ can be inflated to $W$ via $\ell|_{\Omega}=0$.
\end{Pt*}

\subsection{Idempotents and component algebras}

\begin{Pt*}
We have the finite diagonal torus $$\bbT:=\bfT(\bbF_q)$$
and its group ring $R[\bbT]$. As $q-1$ is invertible in $R$, so is $|\bbT|=(q-1)^2$ and hence
$R[\bbT]$ is a semisimple ring. The canonical isomorphism $\bbT\simeq I/I^{(1)}$ induces an inclusion
$$R[\bbT]\subset \cH_{R}^{(1)}(q).$$ We denote by $\bbT^{\vee}$ the set of characters $$\lambda: \bbT\rightarrow \bbF^{\times}_q$$ of $\bbT$, with its natural $W_0$-action given by $$^s\lambda(t_1,t_2)=\lambda(t_2,t_1)$$ for
$(t_1,t_2)\in\bbT$. The number of $W_0$-orbits in $\bbT^{\vee}$ equals $\frac{q^2-q}{2}$.
Also $W$ acts on
$\bbT^{\vee}$ through the canonical quotient map $W\rightarrow W_0$.

\end{Pt*}

\begin{Def*}
For all $\lambda\in \bbT^{\vee}$, define
$$
\varepsilon_{\lambda}:=|\bbT|^{-1}\sum_{t\in\bbT}\lambda^{-1}(t)T_t\quad\in R[\bbT]
$$
and for all $\gamma\in \bbT^{\vee}/W_0$,
$$
\varepsilon_{\gamma}:=\sum_{\lambda\in\gamma}\varepsilon_{\lambda}\quad\in R[\bbT].
$$

\end{Def*}

Following the terminology of \cite{V04}, we call $|\gamma|=1$ the {\it Iwahori case} or {\it non-regular case} and $|\gamma|=2$ the {\it regular case}.

\begin{Prop*}
For all $\lambda\in \bbT^{\vee}$, the element $\varepsilon_{\lambda}$ is an idempotent. For all $\gamma\in \bbT^{\vee}/W_0$, the element $\varepsilon_{\gamma}$ is a central idempotent in $\cH^{(1)}_R(q)$. The $R$-algebra $\cH^{(1)}_R(q)$ is the direct product of its sub-$R$-algebras $\cH^{(1)}_R(q)\varepsilon_{\gamma}$, i.e.
$$
\cH^{(1)}_R(q)=\prod_{\gamma\in\bbT^{\vee}/W_0}\cH^{(1)}_R(q)\varepsilon_{\gamma}.
$$

\end{Prop*}
\begin{proof}
This follows from \cite[Prop. 3.1]{V04} and its proof.
\end{proof}
The proposition implies that the category of $\cH_R^{(1)}(q)$-modules decomposes into a finite product of the module categories for the individual component rings $\cH^{(1)}_R(q)\varepsilon_{\gamma}$.

\subsection{The Iwahori-Hecke algebra}
Our reference for the following is \cite[1.1/2]{V04}. 

\begin{Def*} \label{genericIHalgebra}
Let $\bfq$ be an indeterminate. The \emph{generic Iwahori-Hecke algebra} is the $\bbZ[\bfq]$-algebra $\cH(\bfq)$ defined by generators
$$
\cH(\bfq):=\bigoplus_{(n_1,n_2)\in\bbZ^2}\ \bbZ[\bfq] T_{e^{(n_1,n_2)}}\oplus\bbZ[\bfq] T_{e^{(n_1,n_2)}s}
$$
and relations:

\begin{itemize}
\item braid relations
$$
T_wT_{w'}=T_{ww'}\quad\textrm{for $w,w'\in W$ if $\ell(w)+\ell(w')=\ell(ww')$}
$$
\item quadratic relations
$$
\left\{ \begin{array}{l}
T_s^2=(\bfq-1)T_s+\bfq \\
T_{s_0}^2=(\bfq-1)T_{s_0}+\bfq.
\end{array} \right.
$$

\end{itemize}
\end{Def*}

\begin{Pt*} \label{explicitIwahoriGL2}
Setting $S:=T_s$ and $U:=T_u$, one can check that
$$
\cH(\bfq)=\bbZ[\bfq][S,U^{\pm 1}],\quad S^2=(\bfq-1)S+\bfq,\quad U^2S=SU^2
$$
is a presentation of $\cH(\bfq)$.
For example, $S_0:= T_{s_0}=USU^{-1}$. We also have the generic finite and affine Hecke algebras $$\cH_{0}(\bfq)=\bbZ[\bfq][S]\subset \cH_{\aff}(\bfq)=\bbZ[\bfq][S_0,S].$$ The algebra $\cH_{0}(\bfq)$ has two characters corresponding to $S\mapsto 0$ and
$S\mapsto -1$. Similarly, $\cH_{\aff}(\bfq)$ has four characters. The two characters different from the trivial character $S_0,S\mapsto 0$ and the sign character $S_0,S\mapsto -1$ are called {\it supersingular}.

\end{Pt*}
\begin{Pt*}
The center $Z(\cH(\bfq))$ of the algebra $\cH(\bfq)$ admits the explicit description via
the algebra isomomorphism
\begin{eqnarray*}
\sZ(\bfq):\bbZ[\bfq]\big{[} \Lambda^{+}\big{]}= \bbZ[\bfq]\big{[} e^{(1,0)},(e^{(1,1)})^{\pm 1}\big{]}&  \stackrel{\simeq} {\longrightarrow} & Z(\cH(\bfq)) \\
e^{(1,0)} & \lmapsto & \zeta_1:=U(S-(\bfq-1))+SU\\
e^{(1,1)} & \lmapsto & \zeta_2:=U^2.
\end{eqnarray*}
In particular,
$$
Z(\cH(\bfq))=\bbZ[\bfq][US+(1-\bfq)U+SU,U^{\pm 2}]\subset \bbZ[\bfq][S,U^{\pm1}]=\cH(\bfq).
$$
\end{Pt*}

\begin{Pt*}\label{IsoIw}
Now let $\gamma\in\bbT^{\vee}/W_0$ such that $|\gamma|=1$, say $\gamma=\{\lambda\}$. The ring homomorphism $\bbZ[\bfq]\rightarrow R$, $\bfq\mapsto q$, induces an isomorphism of $R$-algebras
$$
 \cH(\bfq)\otimes_{\bbZ[\bfq]}R\stackrel{\simeq} {\longrightarrow} \cH_R^{(1)}(q)\varepsilon_\gamma,\quad T_w\mapsto \varepsilon_\lambda T_w.
$$
\end{Pt*}

\subsection{The second Iwahori-Hecke algebra} \label{genericgammalagebrareg}
Our reference for the following is \cite[2.2]{V04}, as well as \cite{KK86} for the basic theory of the nil Hecke algebra. We keep the notation introduced above.
\begin{Def*}  \label{genericnilHalgebra}
The \emph{generic nil Hecke algebra} is the $\bbZ[\bfq]$-algebra $\cH^{\nil}(\bfq)$ defined by generators
$$
\cH^{\nil}(\bfq):=\bigoplus_{(n_1,n_2)\in\bbZ^2}\ \bbZ[\bfq] T_{e^{(n_1,n_2)}}\oplus\bbZ[\bfq] T_{e^{(n_1,n_2)}s}
$$
and relations:

\begin{itemize}
\item braid relations
$$
T_wT_{w'}=T_{ww'}\quad\textrm{for $w,w'\in W$  if $\ell(w)+\ell(w')=\ell(ww')$}
$$
\item quadratic relations
$$
\left\{ \begin{array}{l}
T_s^2=\bfq \\
T_{s_0}^2=\bfq.
\end{array} \right.
$$
\end{itemize}
\end{Def*}

\begin{Pt*} \label{explicitnilGL2}
Setting $S:=T_s$ and $U:=T_u$, one can check that
$$
\cH^{\rm nil}(\bfq)=\bbZ[\bfq][S,U^{\pm 1}],\quad S^2=\bfq,\quad U^2S=SU^2
$$
is a presentation of $\cH^{\nil}(\bfq)$. Again, $S_0:=T_{s_0}=USU^{-1}$. The center $Z(\cH^{\rm nil}(\bfq))$ admits the explicit description via the algebra isomorphism
\begin{eqnarray*}
\sZ^{\nil}(\bfq):\bbZ[\bfq]\big{[} \Lambda^{+}\big{]}= \bbZ[\bfq]\big{[} e^{(1,0)},(e^{(1,1)})^{\pm 1}\big{]}& \stackrel{\simeq} {\longrightarrow} & Z(\cH^{\rm nil}(\bfq)) \\
e^{(1,0)} & \lmapsto & \zeta_1:=US+SU\\
e^{(1,1)} & \lmapsto & \zeta_2:=U^2.
\end{eqnarray*}
In particular,
$$
Z(\cH^{\rm nil}(\bfq))=\bbZ[\bfq][US+SU,U^{\pm 2}]\subset \bbZ[\bfq][S,U^{\pm1}]=\cH^{\rm nil}(\bfq).
$$
\end{Pt*}

\begin{Pt*}\label{twistedring} 
Form the twisted tensor product algebra
$$
\cH_2(\bfq):=(\bbZ[\bfq]\times\bbZ[\bfq])\otimes'_{\bbZ[\bfq]} \cH^{\nil}(\bfq).
$$
With the formal symbols $\varepsilon_{1}=(1,0)$ and $\varepsilon_{2}=(0,1)$, the ring multiplication is given by
 $$
 (\varepsilon_i \otimes T_w)\cdot (\varepsilon_{i'} \otimes T_{w'}) = (\varepsilon_i\varepsilon_{{}^wi'}\otimes T_wT_{w'})
 $$
for all $1\leq i,i'\leq 2$. Here, $W$ acts through its quotient $W_0$ and $s\in W_0$ acts on the set $\{1,2\}$ by interchanging the two elements. The multiplicative unit element in the ring $\bbZ[\bfq]\times\bbZ[\bfq]$ is $(1,1)= \varepsilon_{1}+\varepsilon_{2}$ and the multiplicative unit element in the ring $\cH_2(\bfq)$ is $(1,1)\otimes 1$. We identify the rings $\bbZ[\bfq]\times\bbZ[\bfq]$ and $\cH^{\nil}(\bfq)$ with subrings of $\cH_2(\bfq)$ via the maps $(a,b)\mapsto (a,b)\otimes 1$ and $a\mapsto (1,1)\otimes a$ respectively.
In particular, we will write $\varepsilon_1, \varepsilon_2, S_0, S, U\in \cH_2(\bfq)$ etc.

\vskip5pt

We also introduce the generic affine Hecke algebra $$ \cH_{2,\aff}(\bfq)=(\bbZ[\bfq]\times\bbZ[\bfq])\otimes'_{\bbZ[\bfq]}\bbZ[\bfq][S_0,S].$$
It is a subalgebra of $\cH_2(\bfq)$ and has two {\it supersingular} characters $\chi_{1}$ and $\chi_{2}$, namely
$\chi_{1}(\varepsilon_{1})=1$ and $\chi_{1}(\varepsilon_{2})=0$ and
  $\chi_{1}(S_0)=\chi_{1}(S)=0$.
 Similarly for $\chi_{2}$.
\end{Pt*}
\begin{Pt*}\label{matrix} The structure of $\cH_2(\bfq)$ as an algebra over its center can be made explicit.
In fact,
there is an algebra isomorphism with an algebra of $2\times2$-matrices
$$
\cH_2(\bfq)\simeq M(2,\cZ(\bfq)),\quad \cZ(\bfq):=\bbZ[\bfq][X,Y,z_2^{\pm1}]/(XY)
$$
which maps the center $Z(\cH_2(\bfq))$ to the scalar matrices $\cZ(\bfq)$. Under this isomorphism, we have
$$
S\mapsto \left (\begin{array}{cc}
0 & Y \\
z_2^{-1}X & 0
\end{array} \right),
\quad
U\mapsto\left (\begin{array}{cc}
0 & z_2 \\
1 & 0
\end{array} \right),
$$
$$
\varepsilon_{1}\mapsto\left (\begin{array}{cc}
1 & 0 \\
0 & 0
\end{array} \right),
\quad
\varepsilon_{2}\mapsto\left (\begin{array}{cc}
0 & 0 \\
0 & 1
\end{array} \right).
$$
The induced map $Z(\cH_2(\bfq))\ra \cZ(\bfq)$ satisfies
$$
\zeta_1\mapsto\left (\begin{array}{cc}
X+Y & 0 \\
0 & X+Y
\end{array} \right),
\quad
\zeta_2\mapsto\left (\begin{array}{cc}
z_2 & 0 \\
0 & z_2
\end{array} \right).
$$
In particular, the subring
$$
 Z^{\circ}(\cH_2(\bfq)):= \bbZ[\bfq][\zeta_1,\zeta_2^{\pm1}] = Z(\cH^{\nil}(\bfq))\subset  \cH^{\rm nil}(\bfq)\subset\cH_2(\bfq)
$$
lies in fact in the center $Z(\cH_2(\bfq))$ of $\cH_2(\bfq)$.

\end{Pt*}

\begin{Pt*}\label{IsoReg}
Now let $\gamma\in\bbT^{\vee}/W_0$ such that $|\gamma|=2$, say $\gamma=\{\lambda,^s\lambda\}$. The ring homomorphism $\bbZ[\bfq]\rightarrow R$, $\bfq\mapsto q$, induces an isomorphism of $R$-algebras

 $$
\cH_2(\bfq)\otimes_{\bbZ[\bfq]}R \stackrel{\simeq} {\longrightarrow} \cH_R^{(1)}(q)\varepsilon_{\gamma},\quad \varepsilon_1\otimes T_w \mapsto \varepsilon_\lambda T_w,\quad \varepsilon_{2}\otimes T_w \mapsto \varepsilon_{{}^s\lambda} T_w.
 $$
\end{Pt*}

\begin{Rem*} We have used the same letters $S_0,S,U,\zeta_1,\zeta_2$ for the corresponding Hecke operators in the Iwahori Hecke algebra and in the second Iwahori Hecke algebra. This should not lead to confusion, as we will always treat non-regular components and regular components separately in our discussion.
\end{Rem*}

\section{The non-regular case and dual equivariant $K$-theory} \label{IwahoriKtheory}

\subsection{Recollections from algebraic $K^{\widehat{\bfG}}$-theory} \label{recollK}
For basic notions from equivariant algebraic $K$-theory we refer to \cite{Th87}. A useful introduction may also be found in \cite[chap. 5]{CG97}.

\begin{Pt*}

We let  $$\widehat{\bfG}:= {\rm GL}_{2/ \overline{\bbF}_q}$$ be the Langlands dual group of $\bfG$ over the algebraic closure $\overline{\bbF}_q$ of $\bbF_q$. The dual torus
 $$\widehat{\bfT}:=\Spec\; \overline{\bbF}_q[\Lambda]\subset \widehat{\bfG}$$ identifies with the torus of diagonal matrices in $\widehat{\bfG}.$
A basic object is $$R(\widehat{\bfG})\;:= \text{the representation ring of } \widehat{\bfG},$$ i.e. the Grothendieck ring of the abelian tensor category of all finite dimensional $\widehat{\bfG}$-representations. It can be viewed as the equivariant $K$-theory $K^{\widehat{\bfG}}(\pt)$ of the base point
$\pt=\Spec\; \overline{\bbF}_q$. To compute it, we introduce the representation ring $R(\widehat{\bfT})$ of $\widehat{\bfT}$ which identifies canonically, as a ring with $W_0$-action, with the group ring of $\Lambda$, i.e.
$$R(\widehat{\bfT})= \bbZ[\Lambda].$$
The formal character $\chi_V\in \bbZ[\Lambda]^{s}$ of a representation $V$ is an invariant function and is defined by
$$ \chi_V( e^{\lambda} ) = \dim_{\overline{\bbF}_q} V_{\lambda}$$
for all $\lambda\in\Lambda$ where $V_{\lambda}$ is the $\lambda$-weight space of $V$. The map $V\mapsto  \chi_V$ induces a ring isomorphism
$$
 \chi_{\bullet}: R(\widehat{\bfG}) \stackrel{\simeq} {\longrightarrow}\bbZ[\Lambda]^s.
$$
The $\bbZ[\Lambda]^s$-module $\bbZ[\Lambda]$ is free of rank $2$, with basis $\{1,e^{(-1,0)}\}$,
 $$\bbZ[\Lambda]= \bbZ[\Lambda]^s  \oplus \bbZ[\Lambda]^s  e^{(-1,0)}.$$

\end{Pt*}

\begin{Pt*}
We let $$\widehat{\cB}:={\bbP}^1_{\overline{\bbF}_q}$$ be the projective line over $\overline{\bbF}_q$ endowed with its left $\widehat{\bfG}$-action by fractional transformations

$$\left (\begin{array}{cc}
a& b \\
c& d
\end{array} \right)(x) =\frac{ax+b}{cx+d}.$$
Here, $x$ is a local coordinate on ${\bbP}^1_{\overline{\bbF}_q}$. The stabilizer of the point $x=\infty$ is the Borel subgroup $\widehat{\bfB}$ of upper triangular matrices and we may thus write $\widehat{\cB}= \widehat{\bfG}/\widehat{\bfB}$.
We denote by $$K^{\widehat{\bfG}}(\widehat{\cB})\;:=\text{ the Grothendieck group of all } \widehat{\bfG}\text{-equivariant
 coherent }\cO_{\widehat{\cB}}\text{-modules}.$$ Given a representation $V$ and an equivariant coherent sheaf $\cF$, the diagonal action of $\widehat{\bfG}$ makes $\cF\otimes_{\overline{\bbF}_q} V$ an equivariant coherent sheaf. In this way,  $K^{\widehat{\bfG}}(\widehat{\cB})$ becomes a module over the ring $R(\widehat{\bfG})$.

\vskip5pt
The \emph{characteristic homomorphism} in algebraic $K^{\widehat{\bfG}}$-theory is a ring isomorphism
$$
c^{\widehat{\bfG}}: \bbZ[\Lambda] \stackrel{\simeq} {\longrightarrow} K^{\widehat{\bfG}}(\widehat{\cB}).
$$
It maps $e^{\lambda}$ with $\lambda=(\lambda_1,\lambda_2) \in\Lambda$ to the class of the $\widehat{\bfG}$-equivariant line bundle
$\cO_{{\bbP}^1}(\lambda_1-\lambda_2)\otimes \det ^{\lambda_2}$ where $\det$ is the determinant character of $\widehat{\bfG}$. The characteristic homomorphism is compatible with the character morphism  $\chi_{\bullet}$ , i.e. $c^{\widehat{\bfG}}$ is $\bbZ[\Lambda]^{s}\simeq R(\widehat{\bfG})$-linear.
\end{Pt*}

\begin{Pt*}\label{classicalDemazure}
For the definition of the classical Demazure operators on algebraic $K$-theory we refer to \cite{D73,D74}. The Demazure operators
$$
D_s,D_s'\in\End_{R(\widehat{\bfT})^{s}}(R(\widehat{\bfT}))
$$
are defined by:
$$
D_{s}(a)=\frac{a-s(a)}{1-e^{(1,-1)}}\hskip15pt \text{and} \hskip15pt  D_s'(a)=\frac{a-s(a)e^{(1,-1)}}{1-e^{(1,-1)}}
$$
for $a\in R(\widehat{\bfT})$. They are the projectors on $R(\widehat{\bfT})^{s}e^{(-1,0)}$ along $R(\widehat{\bfT})^{s}$, and on $R(\widehat{\bfT})^{s}$ along $R(\widehat{\bfT})^se^{(1,0)}$, respectively.
In particular $ D_s^2=D_s$ and $D'^2_s=D'_s$. One sets
$$
D_s(\bfq):=D_s-\bfq D_s'\in\End_{R(\widehat{\bfT})^{s}[\bfq]}(R(\widehat{\bfT})[\bfq])
$$
and checks by direct calculation that $$D_s(\bfq)^2=\bfq-(\bfq-1)D_s(\bfq).$$ In particular, we obtain a well-defined $\bbZ[\bfq]$-algebra homomorphism
$$
\sA_0(\bfq):\cH_0(\bfq)=\bbZ[\bfq][S] \lra \End_{R(\widehat{\bfT})^{s}[\bfq]}(R(\widehat{\bfT})[\bfq]),\quad S\lmapsto -D_s(\bfq)
$$
which we call the {\it Demazure representation}.
\end{Pt*}

\subsection{The morphism from $R(\widehat{\bfG})[\bfq]$ to the center of $\cH(\bfq)$} \label{def-strI}

In the following we identify the rings
 $$R(\widehat{\bfG})[\bfq]\simeq \bbZ[\bfq][\Lambda]^s=\bbZ[\bfq][\xi_1,\xi_2^{\pm1}]$$
 via the character isomorphism $\chi_{\bullet}$. We have the $\bbZ[\bfq]$-algebra isomorphism coming via base change from the isomorphism $\xi^+$, cf. \ref{cochar}:
\begin{eqnarray*}
\xi^+:\bbZ[\bfq]\big{[} e^{(1,0)},(e^{(1,1)})^{\pm 1}\big{]} & \stackrel{\simeq} {\longrightarrow}& \bbZ[\bfq][\xi_1,\xi_2^{\pm1}] \\
e^{(1,0)} & \lmapsto & \xi_1\\
e^{(1,1)} & \lmapsto & \xi_2.
\end{eqnarray*}
On the other hand, the source of $\xi^+$ is isomorphic to the center $Z(\cH(\bfq))$ of $\cH(\bfq)$ via the isomorphism $\sZ(\bfq)$, cf. \ref{explicitIwahoriGL2}.
The composition

\begin{eqnarray*}
\sZ(\bfq)\circ (\xi^+)^{-1}: R(\widehat{\bfG})[\bfq] &\stackrel{\simeq} {\longrightarrow} &Z(\cH(\bfq)) \\
\xi_1& \lmapsto &\zeta_1=U(S-(\bfq-1))+SU\\
 \xi_2 & \lmapsto &\zeta_2=U^2
\end{eqnarray*}
is then a ring isomorphism.

\subsection{The extended Demazure representation $\sA(\bfq)$} \label{HIqGL2}

Recall the Demazure representation
$\sA_0(\bfq)$ of the finite algebra $\cH_0(\bfq)$ by $R(\widehat{\bfG})[\bfq]$-linear operators on the $K$-theory
$K^{\widehat{\bfG}}(\widehat{\cB})$, cf. \ref{classicalDemazure}.
We have the following first main result.
\begin{Th*}
There is a unique
ring homomorphism
$$
\xymatrix{
\sA(\bfq):\cH(\bfq) \ar[r] & \End_{R(\widehat{\bfG})[\bfq]}(K^{\widehat{\bfG}}(\widehat{\cB})[\bfq])
}
$$
which extends the ring homomorphism $\sA_0(\bfq)$ and coincides
 on $Z(\cH(\bfq))$ with the isomorphism
\begin{eqnarray*}
Z(\cH(\bfq)) &\stackrel{\simeq} {\longrightarrow} & R(\widehat{\bfG})[\bfq] \\
\zeta_1 & \lmapsto &\xi_1\\
 \zeta_2 & \lmapsto &\xi_2.
\end{eqnarray*}
The homomorphism $\sA(\bfq)$ is injective.
\end{Th*}
\begin{Proof} Such an extension exists if and only if there exists
$$\sA(\bfq)(U)\in\End_{R(\widehat{\bfG})[\bfq]}(K^{\widehat{\bfG}}(\widehat{\cB})[\bfq])$$
satisfying
\begin{enumerate}
\item $\sA(\bfq)(U)$ is invertible ;
\item $\sA(\bfq)(U)^2=\sA(\bfq)(U^2)=\sA(\bfq)(\zeta_2)=\xi_2\Id$ ;
\item \begin{eqnarray*}
\sA(\bfq)(U)\sA_0(\bfq)(S)+(1-\bfq)\sA(\bfq)(U)+\sA_0(\bfq)(S)\sA(\bfq)(U)&=&\sA(\bfq)(US+(1-\bfq)U+SU)\\
&=&\sA(\bfq)(\zeta_1)\\
&=&\xi_1\Id.
\end{eqnarray*}

\end{enumerate}

To find such an operator $\sA(\bfq)(U)$, we write
$$
K^{\widehat{\bfG}}(\widehat{\cB})[\bfq]=R(\widehat{\bfT})[\bfq]=R(\widehat{\bfT})^s[\bfq]\oplus R(\widehat{\bfT})^s[\bfq]e^{(-1,0)},
$$
and use the $R(\widehat{\bfT})^s[\bfq]$-basis $\{1,e^{(-1,0)}\}$ to identify $\End_{R(\widehat{\bfG})[\bfq]}(K^{\widehat{\bfG}}(\widehat{\cB})[\bfq])$ with the algebra of $2\times 2$-matrices over the ring $R(\widehat{\bfT})^s[\bfq]$. Then, by definition,
$$
\sA_0(\bfq)(S)=
\left (\begin{array}{cc}
0& 0 \\
0& -1
\end{array} \right)
+
\bfq
\left (\begin{array}{cc}
1& \xi_1 e^{(-1,-1)} \\
0& 0
\end{array} \right)
=
\left (\begin{array}{cc}
\bfq& \bfq\xi_1 e^{(-1,-1)} \\
0& -1
\end{array} \right).
$$
Hence, if we set
$$
\sA(\bfq)(U)=
\left (\begin{array}{cc}
a& c \\
b& d
\end{array} \right),
$$
we get
$$
\sA(\bfq)(U)^2=e^{(1,1)}\Id \Longleftrightarrow
\left (\begin{array}{cc}
a^2+bc& c(a+d) \\
b(a+d)& d^2+bc
\end{array} \right)
=
\left (\begin{array}{cc}
e^{(1,1)}& 0 \\
0& e^{(1,1)}
\end{array} \right)
$$
and
$$
\sA(\bfq)(U)\sA_0(\bfq)(S)+(1-\bfq)\sA(\bfq)(U)+\sA_0(\bfq)(S)\sA(\bfq)(U)=\xi_1\Id
$$
$$
\Longleftrightarrow
\left (\begin{array}{cc}
(\bfq+1)a+\bfq\xi_1 e^{(-1,-1)} b& \bfq\xi_1 e^{(-1,-1)}(a+d) \\
0& -(\bfq+1)d+\bfq\xi_1 e^{(-1,-1)} b
\end{array} \right)
=
\left (\begin{array}{cc}
\xi_1& 0 \\
0& \xi_1
\end{array} \right).
$$
These two conditions together are in turn equivalent to
$$
\left\{ \begin{array}{lll}
a&=&-d\\
bc&=&e^{(1,1)}-a^2\\
(\bfq+1)a&=&\xi_1-\bfq\xi_1 e^{(-1,-1)}b.
\end{array} \right.
$$
Moreover, in this case, the determinant
$$
ad-bc=-a^2-(e^{(1,1)}-a^2)=-e^{(1,1)}
$$
is invertible. Specialising to $\bfq=0$, we find that there is {\it exactly one} $R(\widehat{\bfG})[\bfq]$-algebra homomorphism
$$
\xymatrix{
\sA(\bfq):\cH(\bfq) \ar[r] & \End_{R(\widehat{\bfG})[\bfq]}(K^{\widehat{\bfG}}(\widehat{\cB})[\bfq]),
}
$$
extending the ring homomorphism $\sA_0(\bfq)$, corresponding to the matrix
$$
\sA(\bfq)(U)=
\left (\begin{array}{cc}
a& c \\
b& d
\end{array} \right)
:=
\left (\begin{array}{cc}
\xi_1& e^{(-1,-1)}\xi_1^2-1 \\
-e^{(1,1)}& -\xi_1
\end{array} \right).
$$

Note that $a,b,c,d\in R(\widehat{\bfT})^s\subset R(\widehat{\bfT})^s[\bfq]$. The injectivity of the map $\sA(\bfq)$ will be proved in the next subsection.
\end{Proof}

\subsection{Faithfulness of $\sA(\bfq)$}\label{injIGL2}
Let us show that the map $\sA(\bfq)$ is injective. It follows from \ref{explicitIwahoriGL2} that the ring $\cH(\bfq)$ is generated by the elements
$$
1,\  S,\  U,\  SU
$$
over its center $Z(\cH(\bfq))=\bbZ[\zeta_1,\zeta_2^{\pm1}][\bfq]$. As the latter is mapped isomorphically to the center $R(\widehat{\bfG})[\bfq]=\bbZ[\xi_1,\xi_2^{\pm1}][\bfq]$ of the matrix algebra $\End_{R(\widehat{\bfG})[\bfq]}(K^{\widehat{\bfG}}(\widehat{\cB})[\bfq])$ by
$\sA(\bfq)$, it suffices to check that the images
$$
1,\ \sA_0(\bfq)(S),\ \sA(\bfq)(U),\ \sA_0(\bfq)(S)\sA(\bfq)(U)
$$
of $1, S, U, SU$ by $\sA(\bfq)$ are free over $R(\widehat{\bfG})[\bfq]$. To ease notation, we will write $\xi$ instead of $\xi_1$ in the following calculation. So let $\alpha,\beta,\gamma,\delta\in R(\widehat{\bfT})^s[\bfq]$ (which is an integral domain) be such that
$$
\alpha
\left (\begin{array}{cc}
1& 0 \\
0& 1
\end{array} \right)
+
\beta
\left (\begin{array}{cc}
\bfq& \bfq\xi e^{(-1,-1)} \\
0& -1
\end{array} \right)
+
\gamma
\left (\begin{array}{cc}
a& c \\
b& -a
\end{array} \right)
+
\delta
\left (\begin{array}{cc}
\bfq& \bfq\xi e^{(-1,-1)} \\
0& -1
\end{array} \right)
\left (\begin{array}{cc}
a& c \\
b& -a
\end{array} \right)
=0.
$$
This is equivalent to the expression
$$
\left (\begin{array}{cc}
\alpha& 0 \\
0& \alpha
\end{array} \right)
+
\left (\begin{array}{cc}
\beta\bfq& \beta\bfq\xi e^{(-1,-1)} \\
0& -\beta
\end{array} \right)
+
\left (\begin{array}{cc}
\gamma a& \gamma c \\
\gamma b& -\gamma a
\end{array} \right)
+
\left (\begin{array}{cc}
\delta(\bfq a+\xi e^{(-1,-1)} b\bfq)& \delta\bfq (c-a\xi e^{(-1,-1)}) \\
-\delta b& \delta a
\end{array} \right)
$$ 
being zero, i.e. to the identity
$$
\left (\begin{array}{cc}
\alpha& 0 \\
0& \alpha
\end{array} \right)
+
\left (\begin{array}{cc}
\beta\bfq& \beta\bfq\xi e^{(-1,-1)} \\
0& -\beta
\end{array} \right)
+
\left (\begin{array}{cc}
\gamma a& \gamma c \\
\gamma b& -\gamma a
\end{array} \right)
+
\left (\begin{array}{cc}
\delta(\xi-a)& \delta\bfq (c-a\xi e^{(-1,-1)}) \\
-\delta b& \delta a
\end{array} \right)
=0.
$$
Then
$$
\left\{ \begin{array}{lll}
\alpha+\beta\bfq+\gamma a+\delta(\xi-a)&=&0\\
(\gamma-\delta)b&=&0\\
\beta\bfq\xi e^{(-1,-1)}+\gamma c+\delta\bfq( c-a\xi e^{(-1,-1)})&=&0 \\
\alpha-\beta+(\delta-\gamma)a&=&0.
\end{array} \right.
$$
As $b\neq 0$, we obtain $\delta=\gamma$ and
$$
\left\{ \begin{array}{lll}
\alpha+\beta\bfq+\gamma \xi&=&0\\
\beta\bfq\xi e^{(-1,-1)}+\gamma ((\bfq+1)c-\bfq \xi e^{(-1,-1)} a)&=&0 \\
\alpha-\beta&=&0.
\end{array} \right.
$$
Hence $\alpha=\beta$ and
$$
\left\{ \begin{array}{lll}
\alpha(\bfq+1)+\gamma \xi&=&0\\
\alpha\bfq\xi e^{(-1,-1)}+\gamma ((\bfq+1)c-\bfq \xi e^{(-1,-1)} a)&=&0. \\
\end{array} \right.
$$
The latter system has determinant
$$
(\bfq+1)((\bfq+1)c-\bfq \xi e^{(-1,-1)} a)-\bfq\xi^2 e^{(-1,-1)},
$$
which is nonzero (its specialization at $\bfq=0$ is equal to $c\neq 0$), whence $\alpha=\gamma=0=\beta=\delta$. This concludes the proof and shows that the map $\sA(\bfq)$ is injective. We record the following two corollaries of the proof.
\begin{Cor*}\label{dimOverCenter} 
The ring $\cH(\bfq)$ is a free $Z(\cH(\bfq))$-module on the basis $1, S, U, SU$.
\end{Cor*}
\begin{Cor*}\label{faithfulatzero} 
The representation $\sA(0)$ is injective.
\end{Cor*}

\subsection{Supersingular modules}
In this section we work at $\bfq=0$ and over the algebraic closure $\overline{\bbF}_q$ of the field $\bbF_q$.

\begin{Pt*}\label{Pt_standard}
Consider the ring homomorphism $\bbZ[\bfq]\rightarrow \overline{\bbF}_q$, $\bfq\mapsto q=0$, and let
$$ \cH_{ \overline{\bbF}_q}=  \cH(\bfq)\otimes_{\bbZ[\bfq]}  \overline{\bbF}_q =   \overline{\bbF}_q [S,U^{\pm 1}].$$

The characters of $\cH_{ \overline{\bbF}_q}$ are parametrised by the set $\{0,-1\}\times \overline{\bbF}^\times_q$ via evaluation on the elements $S$ and $U$.
Let $(\tau_1,\tau_2)\in \overline{\bbF}_q \times  \overline{\bbF}_q^{\times}$.
A {\it standard module} over $\cH_{ \overline{\bbF}_q}$ of dimension $2$ is defined to be a module of type
$$ M_2(\tau_1,\tau_2):=\overline{\bbF}_q m \oplus \overline{\bbF}_q Um,\hskip15pt Sm=-m, \hskip15pt SUm=\tau_1 m, \hskip15pt U^2m=\tau_2 m. $$
The center $Z(\cH_{ \overline{\bbF}_q})=\overline{\bbF}_q[\zeta_1,\zeta_2^{\pm 1}]$ acts on the module $M_2(\tau_1,\tau_2)$ via the character $\zeta_1\mapsto \tau_1, \zeta_2\mapsto \tau_2$. The module $M_2(\tau_1,\tau_2)$ is reducible if and only if $\tau_1^2=\tau_2$. It is called {\it supersingular} if $\tau_1=0$. A supersingular module is thus irreducible. Any simple finite dimensional $\cH_{ \overline{\bbF}_q}$-module is either a character or a standard module \cite[1.4]{V04}.
\end{Pt*}

\begin{Pt*}
Now consider the base change of the representation $\sA:=\sA(0)$ to $\overline{\bbF}_q$
$$
\xymatrix{
\sA_{\overline{\bbF}_q}: \cH_{\overline{\bbF}_q} \ar[r] & \End_{R(\widehat{\bfG})_{\overline{\bbF}_q}}(K^{\widehat{\bfG}}(\widehat{\cB})_{\overline{\bbF}_q})=
\End_{\overline{\bbF}_q[\xi_1,\xi_2^{\pm1}]}(\overline{\bbF}_q[e^{\pm\eta_1},e^{\pm\eta_2}]).
}
$$
Recall that the image of
$
Z(\cH_{\overline{\bbF}_q})=\overline{\bbF}_q[\zeta_1,\zeta_2^{\pm1}]$
is $R(\widehat{\bfG})_{\overline{\bbF}_q}=\overline{\bbF}_q[\xi_1,\xi_2^{\pm1}].$

Let us fix a character $\theta: Z(\cH_{\overline{\bbF}_q})\rightarrow \overline{\bbF}_q$. Following $\cite{V04}$, we call $\theta$ {\it supersingular} if $\theta(\zeta_1)=0$. Consider the base change of $\sA_{\overline{\bbF}_q}$ along $\theta$
$$
\cH_\theta:= \cH_{\overline{\bbF}_q}\otimes_{Z(\cH_{\overline{\bbF}_q})}\overline{\bbF}_q ,\quad K^{\widehat{\bfG}}(\widehat{\cB})_{\theta}:=
K^{\widehat{\bfG}}(\widehat{\cB})_{\overline{\bbF}_q}\otimes_{Z(\cH_{\overline{\bbF}_q})}\overline{\bbF}_q
=K^{\widehat{\bfG}}(\widehat{\cB})_{\overline{\bbF}_q}\otimes_{R(\widehat{\bfG})_{\overline{\bbF}_q}}\bigg(R(\widehat{\bfG})_{\overline{\bbF}_q}
\otimes_{Z(\cH_{\overline{\bbF}_q})}\overline{\bbF}_q\bigg),
$$
$$
\xymatrix{
\sA_{\theta}:\cH_{\theta} \ar[r] & \End_{\overline{\bbF}_q}(K^{\widehat{\bfG}}(\widehat{\cB})_{\theta}).
}
$$
\begin{Prop*} The representation $\sA_{\theta}$ is faithful if and only if $\theta(\zeta_1)^2\neq \theta(\zeta_2)$. In this case, $\sA_{\theta}$ is an algebra isomorphism
$$\sA_{\theta}:\cH_{\theta} \stackrel{\simeq}{\longrightarrow} \End_{\overline{\bbF}_q}(K^{\widehat{\bfG}}(\widehat{\cB})_{\theta}).$$

\end{Prop*}
\begin{Proof}
The discussion in the preceding section \ref{injIGL2} shows that $\cH_{\theta}$ has $\overline{\bbF}_q$-basis given by $1, S, U, SU$.
Moreover, their images
$$
1,\ \sA_\theta (S),\ \sA_\theta (U),\ \sA_\theta(S)\sA_\theta(U)
$$
by $\sA_\theta$ are linearly independent over $\overline{\bbF}_q$ if and only if the scalar
$c= e^{(-1,-1)}\xi_1^2-1\in R(\widehat{\bfG})_{\overline{\bbF}_q}$ does not reduce to zero via $\theta$, i.e. if and only if
$\zeta_{2}^{-1}\zeta_1^2-1\notin\ker\theta$. In this case, the map
$\sA_{\theta}$ is injective and then bijective since $\dim_{\overline{\bbF}_q} K^{\widehat{\bfG}}(\widehat{\cB})_{\theta}=2$.
\end{Proof}
\begin{Cor*}\label{propStandardIw}
The $\cH_{\overline{\bbF}_q}$-module $K^{\widehat{\bfG}}(\widehat{\cB})_{\theta}$ is isomorphic to the standard module $M_2(\tau_1,\tau_2)$ where $\tau_1=\theta(\zeta_1)$ and $\tau_2=\theta(\zeta_2)$. In particular, if $\theta$ is supersingular, then $K^{\widehat{\bfG}}(\widehat{\cB})_{\theta}$ is isomorphic to the unique supersingular $\cH_{\overline{\bbF}_q}$-module with central character $\theta$.
\end{Cor*}
\begin{Proof}
In the case $\tau_1^2\neq \tau_2$, the module $K^{\widehat{\bfG}}(\widehat{\cB})_{\theta}$ is irreducible by the preceding proposition and hence is standard.
In general, it suffices to find $m\in K^{\widehat{\bfG}}(\widehat{\cB})_{\theta}$ with $Sm=-m$ and to verify that $\{m,Um\}$ are linearly independent. For example, $m=e^{\eta_2}$ is a possible choice, cf. below.
\end{Proof}

A "standard basis" for the module $K^{\widehat{\bfG}}(\widehat{\cB})_{\theta}$ comes from the so-called
{\it Pittie-Steinberg basis} \cite{St75} of $\overline{\bbF}_q[e^{\pm\eta_1},e^{\pm\eta_2}]$ over $\overline{\bbF}_q[\xi_1,\xi_2^{\pm1}]$. It is given by
\begin{eqnarray*}
z_e&=&1\\
z_{s}&=&e^{\eta_2}.
\end{eqnarray*}
It induces a basis of $\overline{\bbF}_q[e^{\pm\eta_1},e^{\pm\eta_2}]\otimes_{\overline{\bbF}_q[\xi_1,\xi_2^{\pm1}],\theta}\overline{\bbF}_q$ over $\overline{\bbF}_q$ for any character $\theta$ of $\overline{\bbF}_q[\xi_1,\xi_2^{\pm1}]$. Let $\tau_2=\theta(\xi_2)$. The matrices of $S$, $U$ and $S_0=USU^{-1}$ in the latter basis are
$$
S=
\left (\begin{array}{cc}
0 & 0\\
0 & -1
\end{array} \right),\hskip15pt
U=
\left (\begin{array}{cc}
0 & -\tau_2\\
-1 & 0
\end{array} \right),
\hskip15pt
S_0=
\left (\begin{array}{cc}
-1 & 0\\
0 & 0
\end{array} \right).
$$

The two characters of $\cH_{0,\overline{\bbF}_q}=\overline{\bbF}_q[S]$ corresponding to $S\mapsto 0$ and $S\mapsto -1$
are realized by $z_e$ and $z_s$. From the matrix of $S_0$, we see in fact that the whole affine algebra $\cH_{\aff,\overline{\bbF}_q}:=\overline{\bbF}_q[S_0,S]$ acts on $z_e$ and $z_s$ via the two supersingular characters of $\cH_{\aff,\overline{\bbF}_q}$, cf. \ref{explicitIwahoriGL2}.
\end{Pt*}

\begin{Pt*}
We extend this discussion of the component $\gamma=1$ to any other non-regular component as follows.
Consider the quotient map
$$ \bbT^{\vee}\longrightarrow \bbT^{\vee}/W_0.$$
For any $\gamma\in\bbT^{\vee}/W_0$ define the $\overline{\bbF}_q$-variety
$$
\widehat{\cB}^{\gamma}:=\widehat{\cB}\times\pi^{-1}(\gamma).$$

Suppose $|\gamma|=1$. We have the algebra isomorphism
$
 \cH_{\overline{\bbF}_q}\stackrel{\simeq} {\rightarrow} \cH_{\overline{\bbF}_q}^{(1)}\varepsilon_\gamma
$
from \ref{IsoIw}. It identifies the center $Z( \cH_{\overline{\bbF}_q})$ with the center of $\cH_{\overline{\bbF}_q}^{(1)}\varepsilon_\gamma$.
In this way, we let the component algebra $\cH_{\overline{\bbF}_q}^{(1)}\varepsilon_\gamma$ act on $K^{\widehat{\bfG}}(\widehat{\cB})_{\overline{\bbF}_q}$
and we denote this representation by
 $K^{\widehat{\bfG}}(\widehat{\cB}^{\gamma})_{\overline{\bbF}_q}$. We may then state, in obvious terminology, that any supersingular character $\theta$ of the center of $\cH_{\overline{\bbF}_q}^{(1)}\varepsilon_\gamma$ gives rise to the supersingular irreducible $\cH_{\overline{\bbF}_q}^{(1)}\varepsilon_\gamma$-module
 $K^{\widehat{\bfG}}(\widehat{\cB}^{\gamma})_{\theta}$.

\end{Pt*}

\section{The regular case and dual equivariant intersection theory}

\subsection{Recollections from algebraic $CH^{\widehat{\bfG}}$-theory} \label{recollCH}
For basic notions from equivariant algebraic intersection theory we refer to \cite{EG96} and \cite{Bri97}. As in the case of equivariant $K$-theory, the characteristic homomorphism will make everything explicit. 

\begin{Pt*}
We denote by $\Sym(\Lambda)$ the symmetric algebra of the lattice $\Lambda$ endowed with its natural action of the reflection $s$.
The equivariant intersection theory of the base point $\pt=\Spec\; \overline{\bbF}_q$ canonically identifies
with the ring of invariants
$$
\Sym(\Lambda)^{s}\simeq CH^{\widehat{\bfG}}(\pt),
$$
cf. \cite[sec. 3.2]{EG96}. Recall our basis elements $\eta_1:=(1,0)$ and $\eta_2:=(0,1)$ of $\Lambda$, so that
 $\Sym(\Lambda)=\bbZ[\eta_1,\eta_2]$. We define the invariant elements
$$\xi'_{1}:=\eta_1+\eta_2\hskip15pt\text{and}\hskip15pt\xi'_{2}:=\eta_1\eta_2$$ in
$\Sym(\Lambda)^{s}$. Then $$\Sym(\Lambda)^{s}=\bbZ[\xi'_1,\xi'_2]$$
and, after inverting the prime $2$, the $\Sym(\Lambda)^{s}$-module $\Sym(\Lambda)$ is free of rank $2$, on the basis $\{1,\frac{\eta_1-\eta_2}{2}\}$.
\end{Pt*}
\begin{Pt*}
The \emph{equivariant} \emph{Chern class of line bundles} in the algebraic $CH^{\widehat{\bfG}}$-theory of $\widehat{\cB}$ is a map
$$
c_1^{\widehat{\bfG}}:\Pic^{\widehat{\bfG}}(\widehat{\cB}) \lra CH^{\widehat{\bfG}}(\widehat{\cB})
$$
which is a group homomorphism. Then, the corresponding \emph{characteristic homomorphism} is a ring isomorphism
$$
c^{\widehat{\bfG}}: \Sym(\Lambda)\stackrel{\simeq} {\longrightarrow} CH^{\widehat{\bfG}}(\widehat{\cB}),
$$
which maps $\lambda=(\lambda_1,\lambda_2)\in\Lambda$ to the equivariant Chern class of the line bundle $\cO_{\bbP^1}(\lambda_1-\lambda_2)\otimes\det^{\lambda_2}$ on $\widehat{\cB}=\bbP^1_{\overline{\bbF}_q}$, i.e.
$$
c^{\widehat{\bfG}}(\lambda)=c_1^{\widehat{\bfG}}(\cO_{\bbP^1}(\lambda_1-\lambda_2)\otimes{\det}^{\lambda_2}).
$$
Note here that the algebraic group $\widehat{\bfG}= {\rm GL}_{2/\overline{\bbF}_q}$ is {\it special} (in the sense of \cite[6.3]{EG96}) and the map
$c^{\widehat{\bfG}}$ is therefore already bijective at the integral level \cite[sec. 6.6]{Bri97}. The homomorphism $c^{\widehat{\bfG}}$ is $\Sym(\Lambda)^{s}\simeq CH^{\widehat{\bfG}}(\pt)$-linear.

To emphasize the duality and the analogy with the case of $K$-theory (and to ease notation), we abbreviate from now on
$$S(\widehat{\bfT}):=\Sym(\Lambda)\hskip15pt\text{and}\hskip15pt S(\widehat{\bfG}):=\Sym(\Lambda)^{s}.$$
\end{Pt*}

\begin{Pt*}  \label{DCH}
For the definition of the classical Demazure operators on algebraic intersection theory, we refer to \cite{D73}. The Demazure operators
$$
D_s,D_s'\in\End_{S(\widehat{\bfT})^{s}}(S(\widehat{\bfT}))
$$
are defined by:
$$
 D_s(a)=\frac{a-s(a)}{\eta_1-\eta_2}\hskip15pt\text{and}\hskip15pt D_s'(a)=\frac{a-s(a)(1-(\eta_1-\eta_2))}{\eta_1-\eta_2}
$$
for $a\in S(\widehat{\bfT})$. Then $D_s$ is the projector on $S(\widehat{\bfT})^{s}\frac{\eta_1-\eta_2}{2}$ along $S(\widehat{\bfT})^{s}$, and $(-D_s)+D_s'=s$. In particular, $D_s^2=0$ and $D'^2_s=\id$. One sets
$$
D_s(\bfq):=D_s-\bfq D_s'\in\End_{S(\widehat{\bfT})^{s}[\bfq]}(S(\widehat{\bfT})[\bfq])
$$
and checks by direct calculation that $D_s(\bfq)^2=\bfq^2$. We obtain thus a well-defined $\bbZ$-algebra homomorphism
$$
\sA_0^{\nil}(\bfq):\cH_0^{\nil}(\bfq)=\bbZ[\bfq][S] \lra \End_{S(\widehat{\bfT})^{s}[\bfq]}(S(\widehat{\bfT})[\bfq]),\quad \bfq\lmapsto\bfq^2,\quad S\lmapsto -D_s(\bfq)
$$
which we call the {\it Demazure representation}.
\end{Pt*}

\subsection{The morphism from $S(\widehat{\bfG})[\bfq]$ to the center of $\cH^{\nil}(\bfq)$}\label{def-strnil}
The version of the homomorphism $(\xi^{+})^{-1}$ in the regular case is the $\bbZ[\bfq]$-algebra homomorphism

\begin{eqnarray*}
S(\widehat{\bfG})[\bfq]=\bbZ[\bfq][\xi'_1,\xi'_2] & \lra & \bbZ[\bfq]\big{[} e^{(1,0)},(e^{(1,1)})^{\pm 1}\big{]}\\
\xi'_1 & \lmapsto & e^{(1,0)}  \\
\xi'_2 & \lmapsto & e^{(1,1)}
\end{eqnarray*}
which becomes an isomorphism after inverting $\xi'_2$. Its composition with $\sZ^{\nil}(\bfq)$, cf. \ref{explicitnilGL2}, therefore gives a ring isomorphism

\begin{eqnarray*}
S(\widehat{\bfG})[\bfq][\xi'^{-1}_2] &\stackrel{\simeq} {\longrightarrow} &Z(\cH^{\nil}(\bfq)) \\
\xi'_1& \lmapsto &\zeta_1=US+SU\\
 \xi'_2 & \lmapsto &\zeta_2=U^2.
\end{eqnarray*}

\subsection{The extended Demazure representation $\sA_{\bbF_p}^{\nil}(\bfq)$ at $\bfq=0$} \label{HnilqGL2}
Recall the Demazure representation $\sA_0^{\nil}(\bfq)$ of the finite algebra $\cH_0^{\nil}(\bfq)$
by $S(\widehat{\bfG})[\bfq]$-linear operators on the intersection theory $CH^{\widehat{\bfG}}(\widehat{\cB})$, cf. \ref{DCH}.
In this section we work at $\bfq=0$. We write $\sA_0^{\nil}$ for the specialization of $\sA_0^{\nil}(\bfq)$ at $\bfq=0$.
\vskip5pt
For better readibility we make a slight {\it abuse of notation} and denote the elements $\xi'_i$ by $\xi_i$ in this and the following sections.
Moreover, $p$ will always be an {\it odd} prime.
\begin{Pt*} \label{computeGL2reg}
A ring homomorphism
$$
\xymatrix{
\sA^{\nil}:\cH^{\rm nil} \ar[r] & \End_{S(\widehat{\bfG})}(CH^{\widehat{\bfG}}(\widehat{\cB}))
}
$$
which extends $\sA_{0}^{\nil}$ and which is linear with respect to the above ring homomorphism
$S(\widehat{\bfG})\ra Z(\cH^{\nil})$ does {\it not} exist, even after inverting $\xi_2$. However, there exists a natural good approximation (after inverting the prime $2$). We will explain these points in the following.

\begin{Pt*} An extension of $\sA_{0}^{\nil}$, linear with respect to $S(\widehat{\bfG})\ra Z(\cH^{\nil})$, exists
if and only if there is an operator
$$\sA^{\nil}(U)\in\End_{S(\widehat{\bfG})[{\xi}_2^{-1}]}(CH^{\widehat{\bfG}}(\widehat{\cB})[{\xi}_2^{-1}])$$
 satisfying
\begin{enumerate}
\item $\sA^{\rm nil}(U)$ is invertible ;
\item $\sA^{\rm nil}(U)^2=\sA^{\rm nil}(U^2)=\xi_2\Id$, i.e. $\sA^{\nil}(U)^2=\xi_2\Id$ ;
\item $\sA^{\nil}(U)\sA_{0}^{\nil}(S)+\sA_{0}^{\nil}(S)\sA^{\nil}(U)=\sA^{\nil}(US+SU)= \xi_1\Id.$

\end{enumerate}
Tensoring by $\bbF_p$, we may write
$$
CH^{\widehat{\bfG}}(\widehat{\cB})_{\bbF_p}=
S(\widehat{\bfG})_{\bbF_p}\oplus S(\widehat{\bfG})_{\bbF_p}\frac{\eta_1-\eta_2}{2},
$$
and identify
$\End_{S(\widehat{\bfG})_{\bbF_p}}(CH^{\widehat{\bfG}}(\widehat{\cB})_{\bbF_p})
$
with the algebra of $2\times 2$-matrices over the ring $S(\widehat{\bfG})_{\bbF_p}$. The analogous statements hold after inverting $\xi_2$.
\end{Pt*}

Then, by definition,
$$
\sA_{0,\bbF_p}^{\rm nil}(S)=-D_s=
\left (\begin{array}{cc}
0& -1 \\
0& 0
\end{array} \right).
$$
Hence, if we set
$$
\sA_{\bbF_p}^{\rm nil}(U)=
\left (\begin{array}{cc}
a& c \\
b& d
\end{array} \right),
$$
we obtain
$$
\sA_{\bbF_p}^{\rm nil}(U)^2=\xi_2\Id \Longleftrightarrow
\left (\begin{array}{cc}
a^2+bc& c(a+d) \\
b(a+d)& d^2+bc
\end{array} \right)
=
\left (\begin{array}{cc}
\xi_2& 0 \\
0& \xi_2
\end{array} \right)
$$
and
$$
\sA_{\bbF_p}^{\rm nil}(U)(-D_s)+(-D_s)\sA_{\bbF_p}^{\rm nil}(U)=\xi_1\Id
$$
$$
\Longleftrightarrow
\left\{ \begin{array}{lll}
a&=&-d\\
b&=&-\xi_1,
\end{array} \right.
$$
and then the first system becomes equivalent to the equation
$$
a^2-\xi_1c=\xi_2\in {\bbF_p}[\xi_1,{\xi}_2^{\pm1}].
$$
However, since $\xi_2$ has no square root in the ring $\bbF_p[{\xi}_2^{\pm1}]$, this latter equation has no solution (take $\xi_1=0$ !).
Consequently, there does not exist any matrix
$\sA_{\bbF_p}^{\rm nil}(U)$ with coefficients in $S(\widehat{\bfG})_{\bbF_p}[{\xi}_2^{-1}]$ satisfying conditions 1, 2, 3, above.
\end{Pt*}

As a best approximation, we keep condition $1$ and also condition $3$ (up to sign), but, because of the square root obstruction above, we modify condition $2$ to $\sA_{\bbF_p}^{\rm nil}(\zeta_2)=\xi_2^2$. We can then state our second main result.
\begin{Th*}\label{thm2}
There is a distinguished
ring homomorphism
$$
\xymatrix{
\sA^{\nil}_{\bbF_p}:\cH^{\rm nil}_{\bbF_p} \ar[r] & \End_{S(\widehat{\bfG})_{\bbF_p}[\xi_2^{-1}]}(CH^{\widehat{\bfG}}(\widehat{\cB})_{\bbF_p}[\xi_2^{-1}])
}
$$
which extends the ring homomorphism $\sA^{\nil}_0$ and coincides on  $Z(\cH^{\nil})_{\bbF_p}$ with the homomorphism
\begin{eqnarray*}
Z(\cH^{\nil}_{\bbF_p}) &\stackrel{\simeq} {\longrightarrow} & \bbF_p[\xi_1,\xi^{\pm 2}_2]\subset S(\widehat{\bfG})_{\bbF_p} [\xi_2^{-1}]\\
\zeta_1 & \lmapsto &-\xi_1\\
 \zeta_2 & \lmapsto &\xi_2^2.
\end{eqnarray*}
The homomorphism $\sA^{\nil}_{\bbF_p}$ is injective.
\end{Th*}
\begin{proof}
The discussion preceding the theorem shows that the matrix
$$
\sA_{\bbF_p}^{\nil}(U):=
\left (\begin{array}{cc}
(\frac{\xi_1^2}{2}-\xi_2)& -\xi_1(\frac{\xi_1^2}{4}-\xi_2) \\
\xi_1& -(\frac{\xi_1^2}{2}-\xi_2)
\end{array} \right)
$$
does satisfy the three conditions
\begin{enumerate}
\item $\sA_{\bbF_p}^{\rm nil}(U)$ is invertible ;
\item $\sA_{\bbF_p}^{\rm nil}(U)^2=(\xi_2)^2\Id$ ;
\item $\sA_{\bbF_p}^{\rm nil}(US+SU)=-\xi_1\Id.$
\end{enumerate}
The injectivity part of the theorem will be shown in the next subsection. \end{proof}
\begin{Rem*} The minus sign before $\xi_1$ appearing in the value of $\sA_{\bbF_p}^{\rm nil}$ on $\zeta_1=US+SU$
could be avoided by setting $\sA_0^{\nil}(S):=D_s$ instead of $-D_s$ in the Demazure representation. But we will not do this.
\end{Rem*}

\begin{Rem*} In the Iwahori case, one can check that the action of $U$ coincides with the action of the Weyl element $e^{\eta_1}s$. In the regular case, the action of the element $\eta_1s$ does not satisfy the conditions $1$-$3$ appearing in the above proof. However, the action of $\eta^2_1s$ does and, in fact, its matrix is given by matrix $\sA_{\bbF_p}^{\nil}(U)$. So the choice of the matrix $\sA_{\bbF_p}^{\nil}(U)$ is in close analogy with the Iwahori case. Our chosen extension $\sA_{\bbF_p}^{\nil}$ of $\sA_{0,\bbF_p}^{\nil}$ seems to be distinguished for at least this reason. This observation also shows that the action of $U$ can actually be defined integrally, i.e. before inverting the prime $2$.
\end{Rem*}

\subsection{Faithfulness of $\sA_{\bbF_p}^{\nil}$}\label{injnilGL2}

Let us show that the map $\sA_{\bbF_p}^{\nil}$ is injective. It follows from \ref{explicitnilGL2} that the ring
$\cH_{\bbF_p}^{\nil}$ is generated by the elements
$$
1,\  S,\  U,\  SU
$$
over its center $Z(\cH^{\nil}_{\bbF})=\bbF_p[\zeta_1,\zeta_2^{\pm1}]$. The latter is mapped isomorphically to the subring
$$
\bbF_p[\xi_1,\xi_2^{\pm2}]\subset S(\widehat{\bfG})_{\bbF_p}[\xi_2^{-1}]
$$
of the matrix algebra $\End_{S(\widehat{\bfG})_{\bbF_p}[\xi_2^{-1}]}(CH^{\widehat{\bfG}}(\widehat{\cB})_{\bbF_p})$ by $\sA_{\bbF_p}^{\nil}$. For injectivity, it therefore suffices to show that the images
$$
1,\ \sA_{0,\bbF_p}^{\nil}(S),\ \sA_{\bbF_p}^{\nil}(U),\ \sA_{0,\bbF_p}^{\nil}(S)\sA_{\bbF_p}^{\nil}(U)
$$
of $1, S, U, SU$ under $\sA^{\nil}_{\bbF_p}$ are free over $S(\widehat{\bfG})_{\bbF_p}[\xi_2^{-1}]$. To this end, let $\alpha,\beta,\gamma,\delta\in S(\widehat{\bfG})_{\bbF_p}[\xi_2^{-1}]$ (which is an integral domain) be such that
$$
\alpha
\left (\begin{array}{cc}
1& 0 \\
0& 1
\end{array} \right)
+
\beta
\left (\begin{array}{cc}
0 & -1 \\
0& 0
\end{array} \right)
+
\gamma
\left (\begin{array}{cc}
a& c \\
b& -a
\end{array} \right)
+
\delta
\left (\begin{array}{cc}
0 & -1 \\
0& 0
\end{array} \right)
\left (\begin{array}{cc}
a& c \\
b& -a
\end{array} \right)
=0,
$$
i.e.
$$
\left (\begin{array}{cc}
\alpha& 0 \\
0& \alpha
\end{array} \right)
+
\left (\begin{array}{cc}
0 & -\beta  \\
0& 0
\end{array} \right)
+
\left (\begin{array}{cc}
\gamma a& \gamma c \\
\gamma b& -\gamma a
\end{array} \right)
+
\left (\begin{array}{cc}
-\delta b & \delta a \\
0 & 0
\end{array} \right)
=0.
$$
Then
$$
\left\{ \begin{array}{lll}
\alpha+\gamma a-\delta b&=&0\\
\gamma b&=&0\\
-\beta+\gamma c+\delta a&=&0 \\
\alpha-\gamma a&=&0,
\end{array} \right.
$$
with $\alpha,\beta,\gamma,\delta\in S(\widehat{\bfG})_{\bbF_p}[\xi_2^{-1}]$. Now recall our choice
$$
\sA_k^{\nil}(U)=\left (\begin{array}{cc}
a& c \\
b& -a
\end{array} \right)
:=
\left (\begin{array}{cc}
(\frac{\xi_1^2}{2}-\xi_2)& -\xi_1(\frac{\xi_1^2}{4}-\xi_2) \\
\xi_1& -(\frac{\xi_1^2}{2}-\xi_2)
\end{array} \right).
$$
In particular, $b=\xi_1$ implies $\gamma=0$, and then $\alpha=0$, $\delta=0$ and $\beta=0$.
This shows that the map $\sA^{\nil}$ is injective and concludes the proof. We record the following corollary of the proof.
\begin{Cor*} 
The ring $\cH^{\nil}_{\bbF_p}$ is a free $Z(\cH^{\nil}_{\bbF_p})$-module on the basis $1, S, U, SU$.
\end{Cor*}

\subsection{The twisted representation $\sA_{2,\bbF_p}$}

\begin{Pt*}
In the algebra 
$$
 \cH_{2} := \cH_{2}(0) =(\bbZ\times\bbZ )\otimes'_{\bbZ} \cH^{\nil}
$$
we have the two subrings 
$\cH^{\nil}$ and
$\bbZ\times\bbZ$. The aim of this section is to extend the representation  $\sA_{\bbF_p}^{\nil}$ from $\cH_{\bbF_p}^{\nil}$ to the whole algebra $\cH_{2,\bbF_p}:= \cH_2\otimes_{\bbZ}\bbF_p $.
To this end, we consider the $\overline{\bbF}_q$-variety
$$
\widehat{\cB}^2:=\widehat{\cB}_{1}\coprod\widehat{\cB}_{2},
$$
where $\widehat{\cB}_{1}$ and $\widehat{\cB}_{2}$ are two copies of $\widehat{\cB}$. We have
$$
CH^{\widehat{\bfG}}(\widehat{\cB}^2)=CH^{\widehat{\bfG}}(\widehat{\cB}_{1})\times CH^{\widehat{\bfG}}(\widehat{\cB}_{2}).
$$
After base change to $\bbF_p$, the ring $\cH^{\nil}$ acts $S(\widehat{\bfG})[\xi_2^{-1}]$-linearly on $CH^{\widehat{\bfG}}(\widehat{\cB})[\xi_2^{-1}]$ through the map $\sA_{\bbF_p}^{\nil}$.
We extend this action diagonally to $CH^{\widehat{\bfG}}(\widehat{\cB}^2)[\xi_2^{-1}]$, thus defining a ring homomorphism
$$
\xymatrix{
\diag(\sA^{\nil}_{\bbF_p}):\cH_{\bbF_p}^{\nil} \ar[r] & \End_{S(\widehat{\bfG})_{\bbF_p}[\xi_2^{-1}]}(CH^{\widehat{\bfG}}(\widehat{\cB}^2)_{\bbF_p}[\xi_2^{-1}]).
}
$$
Because of the \emph{twisted} multiplication in the algebra $\cH_{2}$, we need to introduce the permutation action of $W$
$$
\xymatrix{
\perm: W \ar@{->>}[r] & W_0 \ar[r] & \Aut_{S(\widehat{\bfG})}(CH^{\widehat{\bfG}}(\widehat{\cB}^2))
}
$$
which permutes the two factors of $CH^{\widehat{\bfG}}(\widehat{\cB}^2)$.

On the other hand, we can consider the projection $p_{i}$ from $CH^{\widehat{\bfG}}(\widehat{\cB}^2)$ to $CH^{\widehat{\bfG}}(\widehat{\cB}_{i})$ as an $S(\widehat{\bfG})$-linear endomorphism of
$CH^{\widehat{\bfG}}(\widehat{\cB}^2)$, for $i=1,2$. The rule $\varepsilon_{i }\mapsto p_{i}$ defines a ring homomorphism
$$
\xymatrix{
\proj: \bbZ\varepsilon_{1}\times\bbZ\varepsilon_{2}\ar[r] & \End_{S(\widehat{\bfG})}(CH^{\widehat{\bfG}}(\widehat{\cB}^2)).
}
$$
\end{Pt*}

\begin{Prop*}\label{embeddingreg}
There exists a unique ring homomorphism
$$
\xymatrix{
\sA_{2,\bbF_p}:\cH_{2,\bbF_p} \ar[r] & \End_{S(\widehat{\bfG})_{\bbF_p}[\xi_2^{-1}]}(CH^{\widehat{\bfG}}(\widehat{\cB}^2)_{\bbF_p}[\xi_2^{-1}])
}
$$
such that
\begin{itemize}
\item

$\quad \sA_{2,\bbF_p}|_{\cH^{\rm nil}_{\bbF_p}}(T_w)=\diag(\sA^{\rm nil}_{\bbF_p})(T_w)\circ \perm(w)\quad\textrm{for all $w\in W$},$
\item
$\quad\sA_{2,\bbF_p}|_{\bbF_p\varepsilon_{1}\times\bbF_p\varepsilon_{2}}=\proj.$
\end{itemize}
The homomorphism $\sA_{2,\bbF_p}$ is injective.
\end{Prop*}
\begin{proof}
Recall that $W_0$ acts on the set $\{1,2 \}$ by interchanging the two elements and then $W$ acts via its projection to $W_0$.
As $\{\varepsilon_{i}T_w, (i,w)\in\{1,2\} \times W\}$ is a $\bbF_p$-basis of $\cH_{2,\bbF_p}$, such a ring homomorphism is uniquely determined by the formula
$$
\sA_{2,\bbF_p}(\varepsilon_{i}T_w)=p_{i}\circ \diag(\sA_{\bbF_p}^{\rm nil})(T_w) \circ \perm(w).
$$
Conversely, taking this formula as a definition of $\sA_{2,\bbF_p}$, we need to check that the resulting $\bbF_p$-linear map is a ring homomorphism, i.e.
$$
\sA_{2,\bbF_p}((1,1))=\Id
$$
and
$$
\sA_{2,\bbF_p}(\varepsilon_{i}T_w\cdot\varepsilon_{i'}T_{w'})=\sA_{2,\bbF_p}(\varepsilon_{i}T_w)\circ\sA_{2,\bbF_p}(\varepsilon_{i'}T_{w'}).
$$
The first condition is clear because $(1,1)=\varepsilon_{1}+\varepsilon_{2}$ and $p_{i}+p_{^si}=\Id$. Let us check the second condition. If $i'\neq {}^{w^{-1}}i$, i.e. $i\neq {}^{w}i'$, then both sides of the claimed equality vanish. Now assume that $i= {}^{w}i'$. On the left hand side we find
$$
\sA_{2,\bbF_p}(\varepsilon_{i}T_w\cdot\varepsilon_{i'}T_{w'})=\sA_{2,\bbF_p}(\varepsilon_{i}T_wT_{w'}),
$$
while on the right hand side, we find
\begin{eqnarray*}
&&\sA_{2,\bbF_p}(\varepsilon_{i}T_w)\circ\sA_{2,\bbF_p}(\varepsilon_{{}^{w^{-1}}i}T_{w'})\\
&=&p_{i}\circ \diag(\sA_{\bbF_p}^{\rm nil})(T_w) \circ \perm(w)\circ p_{{}^{w^{-1}}i}\circ \diag(\sA_{\bbF_p}^{\rm nil})(T_{w'})\circ\perm(w')\\
&=&p_{i}\circ \diag(\sA_{\bbF_p}^{\rm nil})(T_w) \circ \diag(\sA_{\bbF_p}^{\rm nil})(T_{w'})\circ p_{{}^{(w')^{-1}}({}^{w^{-1}}i)}\\
&=&p_{i}\circ \diag(\sA_{\bbF_p}^{\rm nil})(T_wT_{w'})\circ p_{{}^{(ww')^{-1}}i}.
\end{eqnarray*}
If $\ell(ww')\neq\ell(w)+\ell(w')$, then $T_wT_w'=0$ and both sides vanish. Otherwise $T_wT_w'=T_{ww'}$, so that the left hand side becomes
$$
\sA_{2,\bbF_p}(\varepsilon_{i}T_{ww'})=p_{i}\circ \diag(\sA_{\bbF_p}^{\rm nil})(T_{ww'})\circ\perm(ww'),
$$
and the right hand side
$$
p_{i}\circ \diag(\sA_{\bbF_p}^{\rm nil})(T_{ww'})\circ p_{{}^{(ww')^{-1}}i}.
$$
These two operators are equal. This proves the existence and the uniqueness of the extension $\sA_{2,\bbF_p}$. Its injectivity will be shown in the next subsection.
\end{proof}

\subsection{Faithfulness of $\sA_{2,\bbF_p}$}

Let us show that the map $\sA_{2,\bbF_p}$ is injective. This is equivalent to show that the family
$$
\{\sA_{2,\bbF_p}(\varepsilon_{i}T_w),(i,w)\in\{1,2\}\times W\}
$$
is free over $\bbF_p$. So let $\{n_{i,w}\}\in \bbF_p^{(\{1,2\}\times W)}$ such that
$$
\sum_{i,w}n_{i,w}\sA_{2,\bbF_p}(\varepsilon_{i}T_w)=0.
$$
Let us fix $i_0\in\{1,2\}$. Composing by $p_{i_0}$ on the left, we get
$$
\sum_{w}n_{i_0,w}\sA_{2,\bbF_p}(\varepsilon_{i_0}T_w)=0.
$$
The left hand side can be rewritten as
$$
\sum_{w}n_{i_0,w}p_{i_0}\circ \diag(\sA_{\bbF_p}^{\nil}(T_w))\circ \perm(w)=\sum_{w}n_{i_0,w}p_{i_0}\circ\diag(\sA_{\bbF_p}^{\nil}(T_w))\circ p_{{}^{w^{-1}}i_0}.
$$
Now let us fix $w_0\in W_0$. Composing by $p_{{}^{w_0^{-1}}i_0}$ on the right, we get
$$
\sum_{w\in \Lambda w_0}n_{i_0,w}p_{i_0}\circ\diag(\sA_{\bbF_p}^{\nil}(T_w))\circ p_{{}^{w_0^{-1}}i_0}=0.
$$
Then, for each $w\in \Lambda w_0$, remark that
$$
p_{i_0}\circ\diag(\sA_{\bbF_p}^{\nil}(T_w))\circ p_{{}^{w_0^{-1}}i_0}=\iota_{i_0,{}^{w_0^{-1}}i_0}\circ \sA_{\bbF_p}^{\nil}(T_w)\circ p_{{}^{w_0^{-1}}i_0}
$$
in $\End\big(CH^{\widehat{\bfG}}(\widehat{\cB}_1)[\xi_2^{-1}]\times CH^{\widehat{\bfG}}(\widehat{\cB}_2)[\xi_2^{-1}] \big)$, where 
$\iota_{i_0,{}^{w_0^{-1}}i_0}$ is the canonical map
$$
\xymatrix{
\iota_{i_0,{}^{w_0^{-1}}i_0}:CH^{\widehat{\bfG}}(\widehat{\cB}_{{}^{w_0^{-1}}i_0})[\xi_2^{-1}] \ar@{=}[r] & CH^{\widehat{\bfG}}(\widehat{\cB}_{i_0})[\xi_2^{-1}] \ar@{^{(}->}[r] &CH^{\widehat{\bfG}}(\widehat{\cB}_1)[\xi_2^{-1}]\times CH^{\widehat{\bfG}}(\widehat{\cB}_2)[\xi_2^{-1}]. 
}
$$
As the latter is injective, we get
$$
0=\sum_{w\in \Lambda w_0}n_{i_0,w}\sA_{\bbF_p}^{\nil}(T_w)\circ p_{{}^{w_0^{-1}}i_0}=\sA_{\bbF_p}^{\nil}\big(\sum_{w\in \Lambda w_0}n_{i_0,w}T_w\big)\circ p_{{}^{w_0^{-1}}i_0}.
$$
Finally, as $p_{{}^{w_0^{-1}}i_0}:CH^{\widehat{\bfG}}(\widehat{\cB}^{2})[\xi_2^{-1}]\ra CH^{\widehat{\bfG}}(\widehat{\cB}_{{}^{w_0^{-1}}i_0})[\xi_2^{-1}]$ is surjective, and as $\sA_{\bbF_p}^{\nil}$ is injective, cf. \ref{injnilGL2}, we get $n_{i_0,w}=0$ for all $w\in \Lambda w_0$. This concludes the proof that $\sA_{2,\bbF_p}$ is injective.

\subsection{Supersingular modules}

In this section we work over the algebraic closure $\overline{\bbF}_q$ of the field ${\bbF}_q$.

\begin{Pt*}
Recall from \ref{matrix} that

$$ \cH_{2, \overline{\bbF}_q}=  \cH_{2,\bbF_p}\otimes_{\bbF_p}  \overline{\bbF}_q = (\overline{\bbF}_q\times\overline{\bbF}_q)\otimes'_{\overline{\bbF}_q}
  \overline{\bbF}_q [S,U^{\pm 1}]$$
has the structure of a $2\times 2$-matrix algebra over its center
  $Z(\cH_{2,\overline{\bbF}_q})$.
  Since $\overline{\bbF}_q$ is algebraically closed,
  $Z(\cH_{2,\overline{\bbF}_q})$ acts on any finite-dimensional irreducible
  $\cH_{2,\overline{\bbF}_q}$-module by a character (Schur's lemma).
  Let $\theta$ be a character of $Z(\cH_{2,\overline{\bbF}_q})$.
Then $$\cH_{2,\theta}:= \cH_{2,\overline{\bbF}_q}\otimes_{Z(\cH_{2,\overline{\bbF}_q}),\theta} \overline{\bbF}_q$$
is isomorphic to the matrix algebra $M(2,\overline{\bbF}_q)$. In particular, it is a semisimple (even simple) ring.
\end{Pt*}
\begin{Pt*}
The unique irreducible $\cH_{ 2,\overline{\bbF}_q}$-module with
  central character $\theta$ is called the {\it standard module} with character $\theta$. Its $\overline{\bbF}_q$-dimension is $2$ and it is isomorphic to the standard
  representation $\overline{\bbF}_q^{\oplus 2}$ of the matrix algebra $M(2,\overline{\bbF}_q)$.
  The image of the basis $\{(1,0),(0,1)\} $ of $\overline{\bbF}_q^{\oplus 2}$ is called a {\it standard basis}. A central character $\theta$ is called {\it supersingular} if $\theta(X)=\theta(Y)=0$ (or, equivalently,
  if $\theta(\zeta_1)=0$).
 If $\theta$ is supersingular, then the affine algebra $\cH_{2,\aff, \overline{\bbF}_q}$ acts on the standard basis of the module via the characters $\chi_{1}$ respectively $\chi_{2}$ and the action of $U$ interchanges the two, cf. \ref{twistedring} and \ref{matrix}.

 For more details we refer to \cite[2.3]{V04}.
\end{Pt*}

\begin{Pt*}
Now consider the base change of the representation $\sA_{2,\bbF_p}$ to $\overline{\bbF}_q$
$$
\xymatrix{
\sA_{2,\overline{\bbF}_q}:\cH_{2,\overline{\bbF}_q} \ar[r] & \End_{S(\widehat{\bfG})_{\overline{\bbF}_q}[\xi_2^{-1}]}(CH^{\widehat{\bfG}}(\widehat{\cB}^2)_{\overline{\bbF}_q}[\xi_2^{-1}])
=\End_{\overline{\bbF}_q[\xi_1,\xi_2^{\pm1}]}(\overline{\bbF}_q[\eta_1^{\pm1},\eta_2^{\pm1}]^{\oplus 2}).
}
$$
Recall that the image under the map $\sA_{2,\overline{\bbF}_q}$ of the central subring
$$
Z^{\circ}(\cH_{2,\overline{\bbF}_q})=\overline{\bbF}_q[\zeta_1,\zeta_2^{\pm1}]\subset Z(\cH_{2,\overline{\bbF}_q})
$$
is the subring of scalars
$$
\overline{\bbF}_q[\xi_1,\xi_2^{\pm2}]\subset\overline{\bbF}_q[\xi_1,\xi_2^{\pm1}]=S(\widehat{\bfG})_{\overline{\bbF}_q}[\xi_2^{-1}].
$$
\end{Pt*}

\begin{Pt*}
Let us fix a supersingular central character $\theta$ and denote its restriction to $Z^\circ:=Z^{\circ}(\cH_{2,\overline{\bbF}_q})$ by $\theta$, too. Then consider the $\cH_{2,\overline{\bbF}_q}$-action on the base change
$$
CH^{\widehat{\bfG}}(\widehat{\cB}^2)[\xi_2^{-1}]_{\theta}:=
CH^{\widehat{\bfG}}(\widehat{\cB}^2)_{\overline{\bbF}_q}[\xi_2^{-1}]\otimes_{Z^{\circ}}\overline{\bbF}_q
=CH^{\widehat{\bfG}}(\widehat{\cB}^2)_{\overline{\bbF}_q}[\xi_2^{-1}]
\otimes_{S(\widehat{\bfG})_{\overline{\bbF}_q}[\xi_2^{-1}]}\bigg(S(\widehat{\bfG})_{\overline{\bbF}_q}[\xi_2^{-1}]
\otimes_{Z^{\circ}}\overline{\bbF}_q\bigg).
$$
For the base ring, we have
$$
S(\widehat{\bfG})_{\overline{\bbF}_q}[\xi_2^{-1}]\otimes_{Z^{\circ},\theta}\overline{\bbF}_q=
\overline{\bbF}_q[\xi_1,\xi_2^{\pm1}]\otimes_{\sA_{2,\overline{\bbF}_q},\overline{\bbF}_q[\zeta_1,\zeta_2^{\pm1}],\theta}\overline{\bbF}_q
$$
where $\sA_{2,\overline{\bbF}_q}(\zeta_1)=-\xi_1$ and $\sA_{2,\overline{\bbF}_q}(\zeta_2)=\xi_2^2$.  Now put $\theta(\zeta_2)=:b\in\overline{\bbF}_q^{\times}$. Then
$$
S(\widehat{\bfG})_{\overline{\bbF}_q}[\xi_2^{-1}]\otimes_{Z^{\circ},\theta}\overline{\bbF}_q=\overline{\bbF}_q[\xi_1,\xi_2^{\pm1}]/(\xi_1, \xi_2^2-b)
=\overline{\bbF}_q[\xi_2]/( \xi_2^2-b)=: A 
$$
and so
$$
CH^{\widehat{\bfG}}(\widehat{\cB}^2)[\xi_2^{-1}]_{\theta}=
{\overline{\bbF}_q}[\eta_1^{\pm1},\eta_2^{\pm1}]^{\oplus 2}\otimes_{\overline{\bbF}_q[\xi_1,\xi_2^{\pm1}]}\overline{\bbF}_q[\xi_2]/(\xi_2^2-b) = 
{\overline{\bbF}_q}[\eta_1^{\pm1},\eta_2^{\pm1}]^{\oplus 2}\otimes_{\overline{\bbF}_q[\xi_1,\xi_2^{\pm1}]} A.
$$
Note that the $\overline{\bbF}_q$-algebra $A$ is isomorphic to 
the direct product $\overline{\bbF}_q\times\overline{\bbF}_q$ (the isomorphism depending on the choice of a square root of $b$ in $\overline{\bbF}_q$).
An $A$-basis of $CH^{\widehat{\bfG}}(\widehat{\cB}^2)[\xi_2^{-1}]_{\theta}$ is given by the four elements
$ \{1_{i}, \frac{\eta_1-\eta_2}{2}1_{i}\}_{i=1,2}$
where

$$ 1_{i}\in CH^{\widehat{\bfG}}(\widehat{\cB}_{i})
\subset  CH^{\widehat{\bfG}}(\widehat{\cB}_{1})\times CH^{\widehat{\bfG}}(\widehat{\cB}_{2})= CH^{\widehat{\bfG}}(\widehat{\cB}^2)
$$
is the equivariant Chern class of the structure sheaf on $\widehat{\cB}_{i}$, for $i=1,2$.
The $\overline{\bbF}_q$-dimension of $CH^{\widehat{\bfG}}(\widehat{\cB}^2)[\xi_2^{-1}]_{\theta}$ is therefore $8$ and 
 $\cH_{2,\overline{\bbF}_q}$ acts $A$-linearly. The length of the 
  $\cH_{2,\overline{\bbF}_q}$-module $CH^{\widehat{\bfG}}(\widehat{\cB}^2)[\xi_2^{-1}]_{\theta}$ is $4$ and the central character of any irreducible subquotient is necessarily equal to $\theta$, since this is true by construction after restriction to $Z^{\circ}$. In the following, we compute explicitly a composition series.
\end{Pt*}

\begin{Prop*}
The algebra $\cH_{2,\aff,\overline{\bbF}_q}$ acts on $1_{i}\in CH^{\widehat{\bfG}}(\widehat{\cB}^2)[\xi_2^{-1}]_{\theta}$ by the supersingular character $\chi_{i}$, for $i=1,2$.

\end{Prop*}

\begin{Proof}
The action of $\cH_{2,\overline{\bbF}_q}$ on $CH^{\widehat{\bfG}}(\widehat{\cB}^2)_{\overline{\bbF}_q}[\xi_2^{-1}]$ is defined by the map  $\sA_{2,\overline{\bbF}_q}$. Hence, by definition,
$$
\varepsilon_{i'}\cdot 1_{i}=
\left\{ \begin{array}{ll}
1_{i} & \textrm{if $i'=i$} \\
0 & \textrm{otherwise}.
\end{array} \right.
$$
We calculate $$S\cdot 1_{i}=\diag(-D_s)\circ\perm(s)(1_{i})=\diag(-D_s)1_{^si}=0.$$ Moreover,

$$ U^{-1}\cdot 1_{i}=\diag(U^{-1})\circ\perm(u^{-1})(1_{i})=\diag(U^{-1})1_{^si}=s(\eta_1^{-2}) 1_{^si}=\eta_2^{-2} 1_{^si}$$
and $$D_s(\eta_2^{-2})=\frac{\eta^{-2}_2-\eta^{-2}_1}{\eta_1-\eta_2}=(\eta_1\eta_2)^{-2}\frac{\eta^{2}_1-\eta^{2}_2}{\eta_1-\eta_2}=
(\eta_1\eta_2)^{-2}(\eta_1+\eta_2)=\frac{\xi_1}{\xi_2^{2}}.$$ Therefore,
$$ SU^{-1}\cdot 1_{i}=\diag(-D_s)\circ\perm(s) (\eta_2^{-2}1_{^si})=-D_s(\eta_2^{-2})1_{i}=- \frac{\xi_1}{\xi_2^{2}} 1_{i}=0$$
since $\xi_1=0$ in $CH^{\widehat{\bfG}}(\widehat{\cB}_{i})[\xi_2^{-1}]_{\theta}$. It follows
that $S_0\cdot 1_{i}=USU^{-1}\cdot 1_{i} =0$.
\end{Proof}

\begin{Prop*}
A composition series with simple subquotients of the $\cH_{2,\overline{\bbF}_q}$-module 
$$
CH^{\widehat{\bfG}}(\widehat{\cB}^{2})[\xi_2^{-1}]_{\theta}
$$  
is given by
$$
\{0\}
$$
$$
\subset \overline{\bbF}_q1_{i}\oplus \overline{\bbF}_q(U\cdot 1_{i}) 
$$
$$
\subset A1_{i}\oplus A(U\cdot 1_{i})=A1_{i}\oplus A1_{{}^si}
$$
$$
\subset A1_{i}\oplus A1_{{}^si} \oplus \overline{\bbF}_q(\frac{\eta_1-\eta_2}{2}1_{i})\oplus \overline{\bbF}_q(U\cdot\frac{\eta_1-\eta_2}{2}1_{i})
$$
$$
\subset CH^{\widehat{\bfG}}(\widehat{\cB}^{2})[\xi_2^{-1}]_{\theta}.
$$
Here the direct sums $\oplus$ are taken in the sense of
$\overline{\bbF}_q$-vector spaces.
\end{Prop*}

\begin{proof}
First of all,
$$
U\cdot 1_{i}:=\diag(U)\circ\perm(u)(1_{i})=\diag(U)1_{{}^si}=\eta_1^21_{{}^si}=-\xi_21_{{}^si}\in A^{\times} 1_{{}^si}
$$
because $0=\xi_1=\eta_1+\eta_2$ and $0=\xi_1^2=\eta_1^2+\eta_2^2+2\xi_2$ in $CH^{\widehat{\bfG}}(\widehat{\cB})[\xi_2^{-1}]_{\theta}$. Hence the three first $\oplus$ appearing in the statement of the proposition are indeed \emph{direct} sums. These three sums are $U$-stable by construction. Moreover, by the preceding proposition, 
$\cH_{2,\aff,\overline{\bbF}_q}$ acts by the character $\chi_{i}$ on $1_{i}$, hence by the character $\chi_{{}^si}$ on $U\cdot1_{i}$. It follows that $\overline{\bbF}_q 1_{i}\oplus \overline{\bbF}_q(U\cdot 1_{i})$ realizes the standard $\cH_{2,\overline{\bbF}_q}$-module with central character 
$\theta$, and that $A1_{i}\oplus A(U\cdot 1_{i})$ is an $\cH_{2,\overline{\bbF}_q}$-submodule of $CH^{\widehat{\bfG}}(\widehat{\cB}^2)[\xi_2^{-1}]_{\theta}$ of dimension $4$ over $\overline{\bbF}_q$. In fact, if $L\subset A$ is any $\overline{\bbF}_q$-line, the same arguments show that $ L 1_i \oplus L(U\cdot 1_{i})$ realizes the standard $\cH_{2,\overline{\bbF}_q}$-module with central character $\theta$.
In particular, the module $A1_{i}\oplus A(U\cdot 1_{i})$ is semisimple.

\vskip5pt

Now let us compute the action of $\cH_{2,\overline{\bbF}_q}$ on the element $\frac{\eta_1-\eta_2}{2}1_{i}$, for $i=1,2$. We have
$$
\varepsilon_{i'}\cdot \frac{\eta_1-\eta_2}{2}1_{i}=
\left\{ \begin{array}{ll}
\frac{\eta_1-\eta_2}{2}1_{i} & \textrm{if $i'=i$} \\
0 & \textrm{otherwise}.
\end{array} \right.
$$
Next
$$
S\cdot \frac{\eta_1-\eta_2}{2}1_{i}:=\diag(S)\circ\perm(s)(\frac{\eta_1-\eta_2}{2}1_{i})=\diag(S)(\frac{\eta_1-\eta_2}{2}1_{{}^si})=-1_{{}^si},
$$
$$
U^{-1}\cdot \frac{\eta_1-\eta_2}{2}1_{i}:=\diag(U^{-1})\circ\perm(u^{-1})(\frac{\eta_1-\eta_2}{2}1_{i})=\diag(U^{-1})(\frac{\eta_1-\eta_2}{2}1_{{}^si})=\eta_2^{-2}\frac{\eta_2-\eta_1}{2}1_{{}^si},
$$
$$
D_s(\eta_2^{-2}\frac{\eta_2-\eta_1}{2})=\frac{1}{\eta_1-\eta_2}(\eta_2^{-2}\frac{\eta_2-\eta_1}{2}-\eta_1^{-2}\frac{\eta_1-\eta_2}{2})=-\frac{\xi_1^2-2\xi_2}{2\xi_2^2},
$$
$$
SU^{-1}\cdot \frac{\eta_1-\eta_2}{2}1_{i}=\diag(S)\circ\perm(s)(\eta_2^{-2}\frac{\eta_2-\eta_1}{2}1_{{}^si})=\frac{\xi_1^2-2\xi_2}{2\xi_2^2}1_{i},
$$
$$
S_0\cdot\frac{\eta_1-\eta_2}{2}1_{i}:=USU^{-1}\cdot\frac{\eta_1-\eta_2}{2}1_{i}=\diag(U)\circ\perm(u)(\frac{\xi_1^2-2\xi_2}{2\xi_2^2}1_{i})=\eta_1^2\frac{\xi_1^2-2\xi_2}{2\xi_2^2}1_{{}^si}=1_{{}^si}
$$
because $\xi_1=0$ and (hence) $\eta_1^2=-\xi_2$ in $CH^{\widehat{\bfG}}(\widehat{\cB})[\xi_2^{-1}]_{\theta}$, and finally
$$
U\cdot\frac{\eta_1-\eta_2}{2}1_{i}=\diag(U)\circ\perm(u)(\frac{\eta_1-\eta_2}{2}1_{i})=\diag(U)(\frac{\eta_1-\eta_2}{2}1_{{}^si})
=\xi_2\frac{\eta_1-\eta_2}{2}1_{{}^si}
$$
which lies in $A^{\times}(\frac{\eta_1-\eta_2}{2}1_{{}^si})$.
Neither of the two elements $\frac{\eta_1-\eta_2}{2}1_{i}$ and $U\cdot\frac{\eta_1-\eta_2}{2}1_{i}$ lies in the (semisimple) module $A1_{i}\oplus A(U\cdot 1_{i})$. Hence the three last $\oplus$ appearing in the statement are indeed direct and they form a sub-$\cH_{2,\overline{\bbF}_q}$-module of dimension $6$ over $\overline{\bbF}_q$. So the series appearing in the statement is indeed a composition series with irreducible subquotients.
\end{proof}

\begin{Rem*}
We see from the proof of the preceding proposition that the characters of $\cH_{2,\aff,\overline{\bbF}_q}$ in the sub-$\cH_{2,\overline{\bbF}_q}$-module
$$
A1_{i}\oplus A1_{{}^si} \oplus \overline{\bbF}_q(\frac{\eta_1-\eta_2}{2}1_{i})\oplus \overline{\bbF}_q(U\cdot\frac{\eta_1-\eta_2}{2}1_{i})
$$
are contained in $A1_{i}\oplus A1_{{}^si}$. Hence this submodule is \emph{not} semi-simple. \emph{A fortiori} the whole module $CH^{\widehat{\bfG}}(\widehat{\cB}^{2})[\xi_2^{-1}]_{\theta}$ is not semisimple and, hence, has no central character.
\end{Rem*}

\begin{Pt*}
Now we transfer this discussion to any regular component of the algebra  $\cH_{\overline{\bbF}_q}^{(1)}$ as follows.
Let $\gamma=\{\lambda,{}^s\lambda\}\in\bbT^{\vee}/W_0$ be a regular orbit and form the $\overline{\bbF}_q$-variety
$$
\widehat{\cB}^{\gamma}=\widehat{\cB}\times\pi^{-1}(\gamma)=\widehat{\cB}_{\lambda}\coprod\widehat{\cB}_{{}^s\lambda},
$$
where $\widehat{\cB}_{\lambda}$ and $\widehat{\cB}_{{}^s\lambda}$ are two copies of $\widehat{\cB}$. 
We have the algebra isomorphism
$
 \cH_{2,\overline{\bbF}_q}\stackrel{\simeq} {\rightarrow} \cH_{\overline{\bbF}_q}^{(1)}\varepsilon_\gamma
$
from \ref{IsoReg}.
In this way, the representation  $\sA_{2,\overline{\bbF}_q}$ induces a representation 
$$
\xymatrix{
\sA^{\gamma}_{\overline{\bbF}_q}: \cH_{\overline{\bbF}_q}^{(1)}\varepsilon_\gamma \ar[r] & \End_{S(\widehat{\bfG})_{\overline{\bbF}_q}[\xi_2^{-1}]}(CH^{\widehat{\bfG}}(\widehat{\cB}^{\gamma})_{\overline{\bbF}_q}[\xi_2^{-1}]).
}
$$
We may then state, in obvious terminology, that any supersingular character $\theta$ of the center of $\cH_{\overline{\bbF}_q}^{(1)}\varepsilon_\gamma$ gives rise to the $\cH_{\overline{\bbF}_q}^{(1)}\varepsilon_\gamma$-module
 $CH^{\widehat{\bfG}}(\widehat{\cB}^{\gamma})[\xi_2^{-1}]_{\theta}$ and that the semisimplification of the latter module equals a direct sum of four copies of the unique supersingular $\cH_{\overline{\bbF}_q}^{(1)}\varepsilon_\gamma$-module with central character $\theta$. 
\end{Pt*}

\section{Tame Galois representations and supersingular modules}

Our reference for basic results on tame Galois representations is \cite{V94}.

\begin{Pt}
Let $\varpi\in o_F$ be a uniformizer and let $f$ be the degree of the residue field extension $\bbF_q/\bbF_p$, i.e. $q=p^f$.
Let ${\rm Gal}(\overline{F}/F)$ denote the absolute Galois group of $F$. Let $\cI\subset {\rm Gal}(\overline{F}/F)$ be its inertia subgroup. We fix an element $\varphi\in {\rm Gal}(\overline{F}/F)$ lifting the Frobenius $x\mapsto x^q$ on $ {\rm Gal}(\overline{F}/F)/\cI$.
The unique pro-$p$-Sylow subgroup of $\cI$ is denoted by $\cP$ (the wild inertia subgroup) and the quotient $\cI/\cP$ is pro-cyclic with pro-order prime to $p$.
We choose a lift $v\in\cI $ of a topological generator for $\cI/\cP$.
Let $\cW\subset {\rm Gal}(\overline{F}/F)$ denote the Weil group of $F$. The quotient group $\cW/\cP$ is
topologically generated by (the images of) $\varphi$ and $v$ and the only relation between these two generators is $\varphi v \varphi^{-1}=v^q.$ There is a topological isomorphism
 $$\cW/\cP\simeq \varprojlim \bbF_{p^n}^{\times}$$ where the projective limit is taken with respect to the norm maps $\bbF_{p^{nm}}^{\times}\ra \bbF_{p^n}^{\times}$. We denote by $\omega_n$
  the projection map  $\cW/\cP\ra \bbF_{p^n}^{\times}$ followed by the inclusion $\bbF_{p^n}^{\times}\subseteq \overline{\bbF}_q^{\times}$.
  We shall only be concerned with the characters $\omega_f$ and $\omega_{2f}$. The character $\omega_f$ extends from $\cW$ to ${\rm Gal}(\overline{F}/F)$ by choosing a root $\sqrt[q-1]{-\varpi}$ and letting ${\rm Gal}(\overline{F}/F)$ act as
  $$ g\mapsto\frac{g\sqrt[q-1]{-\varpi}}{\sqrt[q-1]{-\varpi}}\in \mu_{q-1}(F)$$
  followed by reduction mod $\varpi$. The character $$\omega_f: {\rm Gal}(\overline{F}/F)\longrightarrow \bbF_{q}^{\times}$$ depends on the choice of $\varpi$ (but not on the choice of $\sqrt[q-1]{-\varpi}$) and equals the reduction mod $\varpi_F$ of the Lubin-Tate character $\chi_L: {\rm Gal}(\overline{F}/F)\rightarrow o_F^\times$ associated to the uniformizer $\varpi$. By changing $\varphi$ by an element of $\cI$, if necessary, we may assume $\omega_f(\varphi)=1$.
  We normalize local class field theory $\cW^{\rm ab}\simeq F^\times$ by sending the geometric Frobenius $\varphi^{-1}$ to $\varpi$. We view the restriction of $\omega_f$ to $\cW$ as a character of $F^\times$.
\end{Pt}

\begin{Pt}

The set of isomorphism classes of irreducible smooth Galois representations
$$\rho : {\rm Gal}(\overline{F}/F)\longrightarrow \widehat{\bfG}={\rm GL}_2(\overline{\bbF}_q)$$ is in bijection with the set of equivalence classes of pairs $(s,t)\in\widehat{\bfG} ^{2}$ such that
$$ s=
\begin{pmatrix}
0  & 1    \\
-b  & 0

\end{pmatrix} \hskip30pt {\rm and } \hskip30pt
t=
\begin{pmatrix}
y & 0    \\
0  & y^q
\end{pmatrix}
$$
with $b\in\overline{\bbF}_q^{\times}$ and $y\in\bbF_{q^2}\setminus \bbF_{q}.$ Here, two pairs $(s,t)$ and $(s',t')$ are equivalent if $s=s'$ and $t,t'$ are ${\rm Gal}(\bbF_{q^2}/\bbF_{q})$-conjugate. Note that $\det(s)=b$ and that $sts^{-1}=t^q$. The bijection is induced by the map $\rho\mapsto (\rho(\varphi),\rho(v))$. The number of equivalence classes of such pairs $(s,t)$ equals $\frac{q^2-q}{2}$ and hence coincides with the number of $W_0$-orbits in $\bbT^{\vee}$.

\vskip5pt
\end{Pt}
\begin{Pt}
By the above numerical coincidence (the "miracle" from \cite{V04}), there exist (many) bijections between the isomorphism classes of irreducible smooth two-dimensional Galois representations and the isomorphism classes of supersingular two-dimensional $\cH^{(1)}_{\overline{\bbF}_q}$-modules. In the following we discuss a a certain example of such a bijection in our geometric language.

\vskip5pt

Let $\rho$ be a two-dimensional irreducible smooth Galois representation with parameters $(s,t)$. Since the element $\omega_{2f}(v)$ generates $\bbF_{q^2}^{\times}$, the element $t$ uniquely determines an exponent $1\leq h\leq q^2-1$, such that
 $$\omega_{2f}(v)^h=y.$$ Replacing $\rho$ by an isomorphic representation $\rho'$ which replaces $y$ by its Galois conjugate $y^q$ replaces $h$ by the rest of the euclidian division of $qh$ by $q^2-1$. We call either of the two numbers an {\it exponent} of $\rho$. 
 
  \begin{Lem} There is $0\leq i\leq q-2$ such that $\rho\otimes \omega_{f}^{-i}$ has an exponent $\leq q-1$.
 \end{Lem}
 \begin{Proof}
This is implicit in the discussion in \cite{V94}. Let $\omega_{2f}(v)^h=y$. Then $h\leq q^2-2$ since $y\neq 1$.
Moreover, $q^2-2 - (q-2)(q+1) = q$. Since $\omega_{2f}^{q+1}=\omega_f$, twisting with $\omega_f$ reduces to the case $h\leq q$.
 Replacing $y$ by its Galois conjugate $y^q$, if necessary, reduces then further to $h\leq q-1$.
 \end{Proof}

By the lemma, we may associate two numbers $1\leq h\leq q-1$  and $0\leq i\leq q-2$ to the representation $\rho$.
We form the character
$$\omega_f^{h-1+i}\otimes \omega_f^i: (F^{\times})^2\longrightarrow \bbF^{\times}_q,\ (t_1,t_2)\mapsto\omega_f^{h-1+i}(t_1) \omega_f^i(t_2) $$
and restrict to $\mu_{q-1}(F)^2.$ This gives rise to an element $\lambda(\rho)$ of $\bbT^{\vee}$ and we take its $W_0$-orbit $\gamma_\rho$.
 \begin{Lem}
 The orbit $\gamma_\rho$ depends only on the isomorphism class of $\rho$.
 \end{Lem}
 \begin{Proof} Suppose $\rho'\simeq \rho$ with
 $$\rho'(v)=t'=
\begin{pmatrix}
y^q & 0    \\
0  & y
\end{pmatrix}.
$$
By the preceding lemma, there is $0\leq i\leq q-2$
and an exponent $1\leq h\leq q-1$ of $\rho\otimes\omega_f^{-i}$. 
If $1<h$, then by definition $\omega_{2f}^h(v)=y\omega_{f}^{-i}(v)$, so that $\omega_{2f}^{qh}(v)=y^q\omega_{f}^{-i}(v)$, and hence $\omega_{2f}^{q-(h-1)}(v)=y^q\omega_{f}^{-(h-1+i)}(v)$, using  $qh=q-(h-1)+(h-1)(q+1)$. Then $1\leq h':=q-(h-1) \leq q-1$  and taking $0\leq i'\leq q-2$ congruent to $h-1+i~ {\rm mod}~ q-1$, we obtain that $h'$ is an exponent for $\rho'\otimes\omega_f^{-i'}$. In particular, $\lambda(\rho'):=\omega_f^{h'-1+i'}\otimes \omega_f^{i'}$, which is $s$-conjugate to $\lambda(\rho)$. If $h=1$, then by definition $\omega_{2f}(v)=y^q\omega_{f}^{-i}(v)$, which implies $\lambda(\rho')=\lambda(\rho)$ in this case. 
\end{Proof}

We call $\rho$ {\it (non-)regular} if the orbit $\gamma_\rho$ is (non-)regular. On the other hand, we view the element $s=\rho(\varphi)$ as a supersingular character $\theta_\rho$ of the center $Z(\cH_{\overline{\bbF}_q}^{\gamma_\rho})$, i.e. $\theta_{\rho}(\zeta_1)=0$ and $\theta_\rho(\zeta_2)=b$. Finally, we have the $\overline{\bbF}_q$-variety
$$ \widehat{\cB}^{\gamma}=\widehat{\cB}\times\pi^{-1}(\gamma)$$
coming from the quotient map $\bbT^{\vee}\ra \bbT^{\vee}/W_0$. These data give rise to the supersingular $\cH^{(1)}_{\overline{\bbF}_q}$-module
$$
\cM(\rho):=
\left\{ \begin{array}{ll}
K^{\widehat{\bfG}}(\widehat{\cB}^{\gamma})_{\theta_\rho} &\text{\rm if }\rho \text{\rm ~ non-regular}\\
&\\
CH^{\widehat{\bfG}}(\widehat{\cB}^{\gamma})[\xi_2^{-1}]_{\theta_\rho}&\text{\rm if }\rho \text{\rm ~regular}.
\end{array} \right.
$$
Recall that $\cH^{(1)}_{\overline{\bbF}_q}$ acts on $\cM(\rho)$ via the projection onto $\cH^{(1)}_{\overline{\bbF}_q}\varepsilon_{\gamma_\rho}$ followed by the extended
Demazure representation $\sA^{\gamma_\rho}_{\overline{\bbF}_q}$. Recall also that the semisimplification of $\cM(\rho)$ is a direct sum of four copies of the supersingular standard module, if $\rho$ is regular. By abuse of notation, we denote a simple subquotient of $\cM(\rho)$ again by $\cM(\rho)$.

\begin{Prop}\label{bijection}
The map $\rho\mapsto\cM(\rho)$ gives a bijection between the isomorphism classes of two-dimensional irreducible smooth $\overline{\bbF}_q$-representations of ${\rm Gal}(\overline{F}/F)$
and the isomorphism classes of two-dimensional supersingular $\cH^{(1)}_{\overline{\bbF}_q}$-modules.
\end{Prop}
\begin{Proof}
By construction, the restriction of $\omega_{f}^{h-1}$ to $\mu_{q-1}(F)\simeq \bbF_q^{\times}$ is given by the exponentiation
$x\mapsto x^{h-1}$. Given $0\leq i\leq q-2$ and $1\leq h\leq q-1$, and $b\in\overline{\bbF}^{\times}_q$, the parameter
$y:=\omega_{2f}(v)^{h}$ lies in $\bbF_{q^2}\setminus \bbF_{q}$
and the pair $(s,t)$ determines a Galois representation $\rho$ having $h$ comme exponent. Hence, $\rho\otimes\omega_f^{i}$ gives rise to the character $\omega_f^{h-1+i}\otimes \omega_f^i$.
The elements of type $\gamma_\rho$ exhaust therefore all orbits in $\bbT^\vee/W_0$. Since a two-dimensional supersingular $\cH^{(1)}_{\overline{\bbF}_q}$-module is determined by its $\gamma$-component and its central character, the map $\rho\mapsto\cM(\rho)$ is seen to be surjective. It is then bijective, since source and target have the same cardinality.
\end{Proof}

\end{Pt}
\begin{Pt}

Let $F$ be a finite extension of $\bbQ_p$. A distinguished natural bijection between 
irreducible 
two-dimensional ${\rm Gal}(\overline{F}/F)$-representations and supersingular two-dimensional $\cH^{(1)}_{\overline{\bbF}_q}$-modules is established by Breuil \cite{Br03} for $F=\bbQ_p$ (see \cite{Be11} for its relation to the $p$-adic local Langlands correspondence for ${\rm GL}_2(\bbQ_p)$) and by Grosse-Kl\"onne \cite{GK18} for general $F/\bbQ_p$. In this final paragraph we will show that the bijection $\rho\mapsto\cM(\rho)$ from \ref{bijection} coincides in this case with the bijections \cite{Br03} and \cite{GK18}.

\vskip5pt

The case $F=\bbQ_p$ follows directly from the explicit formulae given in \cite[1.3]{Be11}. For the general case, we briefly recall the main construction from \cite{GK18} in the case of standard supersingular modules of dimension $2$.
Let $F_{\phi}$ be the special Lubin-Tate group with Frobenius power series 
$\phi(t)=\varpi t + t^q$. Let $F_{\infty}/F$ be the extension generated by all torsion points of $F_{\phi}$ and let $\Gamma = {\rm Gal}(F_{\infty}/F)$.
We identify in the following $\Gamma \simeq o_F^{\times}$ via the character $\chi_L$.

Let $k/\bbF_q$ be a finite extension and let 
$\cH^{(1)}_k:=\cH^{(1)}(\bfq)\otimes_{\bbZ[\bfq]} k$ via $\bfq\mapsto q=0$.
Let $M$ be a two-dimensional standard supersingular $\cH^{(1)}_k$-module, arising from a supersingular character $\chi: \cH^{(1)}_{{\rm aff },k}\rightarrow k$.
Let $e_0\in M$ such that $\cH^{(1)}_{{\rm aff },k}$ acts on $e_0$ via $\chi$ and put $e_1=T^{-1}_{\omega} e_0$ (where $\omega=u^{-1}$ in our notation).\footnote{For example, if $M$ is an $\cH_{\theta}$-module on which $U^2=\zeta_2$ acts via the scalar $\theta(\zeta_2)=\tau_2$, then $U=U^{-1} \cdot \tau_2$ on $M$ and
$m:=\tau_2^{-1} e_1$ satisfies $U m= T_{\omega} e_1 = e_0$, i.e. $\{m,Um\}$ is a standard basis for $M$ in the sense of \ref{Pt_standard}.} The character $\chi$ determines two numbers 
$0\leq k_0,k_1 \leq q-1$ with $(k_0,k_1)\neq (0,0), (q-1,q-1)$. One considers $M$ a $k[[t]]$-module with $t=0$ on $M$. Let $\Gamma = o_F^{\times}$ 
act on $M$ via $$\gamma (m) = T_{\eta_1(\overline{\gamma})}^{-1} (m)$$ for $\gamma\in o_F^{\times}$ with reduction $\overline{\gamma}\in\bbF_q^{\times}$ and 
$\eta_1(\overline{\gamma})^{-1}=\diag(\overline{\gamma}^{-1},1) \in\bbT$. The 
$k[[t]][\varphi]$-submodule $\nabla(M)$ of $$k[[t]][\varphi,\Gamma ]\otimes_{k[[t]][\Gamma]} M \simeq k[[t]][\varphi]\otimes_{k[[t]]} M$$ is then generated by the two elements 
$h(e_j) = t^{k_j}\varphi \otimes T_{\omega}^{-1} (e_j) + 1 \otimes e_j$ thereby defining the relation between the Frobenius $\varphi$ and the Hecke action of $T_{\omega}$.
Note that in the case of ${\rm GL}_2$, the cocharacter $e^*$ of \cite[2.1]{GK18} is equal to $\eta_1.$

The module $\nabla(M)$ is stable under the $\Gamma$-action and thus the quotient
$$\Delta(M) =  \big( k[[t]][\varphi]\otimes_{k[[t]]} M \big)  / \nabla(M)$$
defines a $k[[t]][\varphi,\Gamma ]$-module. It is torsion standard cyclic with weights $(k_0,k_1)$ in the sense of \cite[1.3]{GK18}, according to 
\cite[Lemma 5.1]{GK18}.
Let $\Delta(M) ^{*} = \Hom_k (\Delta(M),k)$. By a general construction
(which goes back to Colmez and Emerton in the case $F=\bbQ_p$ and $\phi(t)=(1+t)^p-1$, as recalled in \cite[2.6]{Br15}) the $k((t))$-vector space 
$$\Delta(M) ^{*} \otimes_{k[[t]]} k((t))$$
is in a natural way an \' etale Lubin-Tate $(\varphi,\Gamma)$-module of dimension $2$. The correspondence $M\mapsto \Delta(M) ^{*} \otimes_{k[[t]]} k((t))$ extends in fact to a fully faithful functor from 
a suitable category of supersingular $\cH^{(1)}_k$-modules to the category of \' etale $(\varphi,\Gamma)$-modules over $k((t))$. The composite functor to the category of continuous ${\rm Gal}(\overline{F}/F)$-representations over $k$ is denoted by $M\mapsto V(M)$. It induces the aforementioned bijection between irreducible 
two-dimensional ${\rm Gal}(\overline{F}/F)$-representations and supersingular two-dimensional $\cH^{(1)}_{\overline{\bbF}_q}$-modules.
\begin{Prop}
The inverse map to the bijection $M\mapsto V(M)$ 
is given by the map $\rho\mapsto \cM(\rho)$.
\end{Prop}
\begin{Proof}
The correspondence $M\mapsto V(M)$ is 
compatible with the twist by a character of $F^\times$ and local class field theory, such that the determinant corresponds to the central character restricted to $F^\times$. By its very construction, the same is true for the correspondence $\rho\mapsto\cM(\rho)$. It therefore suffices to compare them on irreducible Galois representations having parameters $b=1$ and $i=0$. Let $k=\bbF_{q^2}$ in the following.
Let $\ind (\omega_{2f}^h)$ be the Galois representation with exponent $1\leq h \leq q-1$ and $b=1$ and $i=0$. Let
$D$ be the $(\varphi,\Gamma)$-module associated to $\rho:=\ind (\omega_{2f}^h)$ and let $M$ be a supersingular  $\cH^{(1)}_k$-module such that 
$\Delta(M) ^{*} \otimes_{k[[t]]} k((t)) \simeq D$. According to the main result of \cite{PS3} for $n=2$, the module $D$ has a basis $\{g_0,g_1\}$
such that $$\gamma(g_j)=\overline{f}_{\gamma}(t)^{hq^j/(q+1)}g_j$$ for all $\gamma\in\Gamma$ and $\varphi(g_0)=g_{1}$ and $\varphi(g_1)=-t^{-h(q-1)}g_0$.
Here, $\overline{f}_{\gamma}(t)=\omega_f(t)t/\gamma(t)\in k[[t]]^\times$.
Define the triple $(k_0,k_1,k_2)=(h-1, q-h, h-1)$ and let $i_j := q-1 -k_{2-j}$, so that $i_0=i_2=q-h$ and $i_1=2q-h-1$. 
Define the triple $(h_0,h_1,h_2)=(0, i_1, i_0 +i_1 q)$. Note that $h_2=h(q-1).$
Put $f_j=t^{h_j}g_j$ for $j=0,1$ and let $D^{\sharp}\subset D$
be the $k[[t]]$-submodule generated by $\{f_0,f_1\}$. Let $(D^{\sharp})^{\ast}$ be the $k$-linear dual. Define $e_i' \in (D^{\sharp})^{\ast}$ via
$e'_i(f_j)=\delta_{ij}$ and $e'_i =0$ on $t D^{\sharp}$. Using the explicit formulae for the $\psi$-operator on $k((t))$ as described in \cite[Lemma 1.1]{GK18} one may follow the argument of \cite[Lemma 6.4]{GK16} and show that $D^{\sharp}$ is a $\psi$-stable lattice in $D$ and that $\{e'_0,e'_1\}$ is a $k$-basis of the $t$-torsion part of  $(D^{\sharp})^{\ast}$ satisfying

$$ t^{k_1} \varphi (e'_0) = e'_1 \hskip5pt \text{and} \hskip5pt  t^{k_0} \varphi (e'_1) = - e'_0. $$

But according to \cite[1.15]{GK18} there is only one $\psi$-stable lattice in $\Delta(M) ^{*} \otimes_{k[[t]]} k((t))$, 
namely $\Delta(M)^*$. It follows that $\Delta(M)\simeq (D^{\sharp})^{\ast}$ and so the weights of the torsion standard cyclic $k[[t]][\varphi,\Gamma ]$-module
$\Delta(M)$ are $(k_0,k_1)$. Since $k_0=h-1$, one deduces from \cite[Lemma 4.1/5.1]{GK18} that $\epsilon_1 \equiv h-1 \mod (q-1).$ This means $\lambda \circ \alpha^{\vee}(x)^{-1}=x^{h-1}$
for the character $\lambda\in\bbT^{\vee}$ of $M$. Since $i=0$ and hence $a=0$ (in the notation of \cite[2.2]{GK16}), we arrive therefore at 
$$ \lambda (\diag(x_1,x_2)) =\lambda ( e^*(x_1x_2)\alpha^{\vee}(x_2)^{-1}) =e^*(x_1x_2)^{a} x_2^{h-1}=x_2^{h-1}.$$
Hence the image of $\lambda$ in $\bbT^{\vee}/W_0$ coincides with $\gamma_{\rho}$. This implies $M \simeq \cM(\rho)$, as claimed.  
\end{Proof}
\end{Pt}

\vskip10pt 

\noindent {\small Cédric Pépin, LAGA, Université Paris 13, 99 avenue Jean-Baptiste Clément, 93 430 Villetaneuse, France \newline {\it E-mail address: \url{cpepin@math.univ-paris13.fr}} }

\vskip10pt

\noindent {\small Tobias Schmidt, IRMAR, Universit\'e de Rennes 1, Campus Beaulieu, 35042 Rennes, France \newline {\it E-mail address: \url{tobias.schmidt@univ-rennes1.fr}} }


\begin{thebibliography}{99}

\bibitem[Be10]{Be10} {\sc L. Berger},
{\it On some modular representations of the Borel subgroup 
of ${\rm GL_2}(\bbQ_p)$}, Comp. Math. \textbf{146}(1) (2010), 58-80.

\bibitem[Be11]{Be11} {\sc L. Berger},
{\it La correspondance de Langlands locale $p$-adique pour ${\rm GL_2}(\bbQ_p)$}, Astérisque \textbf{339} (2011), 157-180.

\bibitem[Br03]{Br03} {\sc C. Breuil}, {\it Sur quelques représentations modulaires et $p$-adiques de $GL_2(\bbQ_p)$. I}, Compositio Math. \textbf{138} (2003), 165-188.

\bibitem[Br15]{Br15} {\sc C. Breuil}, {\it Induction parabolique et $(\varphi,\Gamma)$-modules.}, Algebra and Number Theory. \textbf{9} (2015), 2241-2291.

\bibitem[Bri97]{Bri97} {\sc M. Brion}, {\it Equivariant Chow groups for torus actions}, Transformations Groups, Vol. 2, Nr. 3, 1997, pages 225-267.

\bibitem[CG97]{CG97} {\sc N. Chriss and V. Ginzburg},
{\it Representation theory and complex geometry}, Birkhäuser, Boston,1997.

\bibitem[C10]{C10} {\sc P. Colmez}, {\it Représentations de $GL_2(\bbQ_p)$ et $(\varphi,\Gamma)$-modules}, Astérisque \textbf{330} (2010), 281-509.

\bibitem[D73]{D73} {\sc M. Demazure}, {\it Invariants symétriques entiers des groupes de Weyl et torsion}, Invent. Math. 21, pages 287-301, 1973.

\bibitem[D74]{D74} {\sc M. Demazure}, {\it Désingularisation des variétés de Schubert généralisées}, Ann. Sci. \'Ecole Norm. Sup. (4), tome 7 (1974), 53-88. Collection of articles dedicated to Henri Cartan on the occasion of his 70th birthday, I.

\bibitem[EG96]{EG96} {\sc D. Edidin and W. Graham}, {\it Equivariant Intersection Theory}, Invent. Math. 131, pages 595-634, 1996.

\bibitem[GK16]{GK16} {\sc E. Grosse-Klönne}, {\it From pro-$p$-Iwahori-Hecke modules to $(\varphi,\Gamma)$-modules, I}, Duke Math. Journal 165 No. 8 (2016), 1529-1595.

\bibitem[GK18]{GK18} {\sc E. Grosse-Klönne}, {\it Supersingular Hecke modules as Galois representations}, Preprint (2018) arXiv:1803.02616.

\bibitem[Hu05]{Hu05} {\sc J.E. Humphreys}, {\it Modular Representations of Finite Groups of Lie Type}, London Math. Soc. Lecture Notes Series 326, Cambridge University Press, 2005.

\bibitem[KK86]{KK86} {\sc B. Kostant, S. Kumar}, {\it The Nil Hecke ring and cohomology of $G/P$ for a Kac-Moody group $G$},
Advances. Math., \textbf{62} (1), 187-237, 1986.

\bibitem[KL87]{KL87} {\sc D. Kazhdan, G. Lusztig}, {\it Proof of the Deligne-Langlands conjecture for Hecke algebras},
Invent. Math., \textbf{87} (1), 153-215, 1987.


\bibitem[KR09]{KR09} {\sc M. Kisin, W. Ren}, {\it Galois representations and Lubin-Tate groups},
Doc. Math., \textbf{14}, 441-461, 2009.

\bibitem[O14]{O14} {\sc R. Ollivier}, {\it Compatibility between Satake and Bernstein isomorphisms in characteristic $p$}, Algebra and Number Theory \textbf{8}(5) (2014), 1071-1111.

\bibitem[PS2]{PS2} {\sc C. Pépin, T. Schmidt}, {\it Mod $p$ Hecke algebras and dual equivariant cohomology II: the case of $GL_n$}, in preparation.

\bibitem[PS3]{PS3} {\sc C. Pépin, T. Schmidt}, {\it On a certain class of Lubin-Tate {$(\varphi,\Gamma)$}-modules}, Preprint 2019.

\bibitem[Sch17]{Sch17} {\sc P. Schneider}, {\it Galois representations and {$(\varphi,\Gamma)$}-modules}, Cambridge Studies in Advanced Mathematics \textbf{164} (Cambridge University Press, 2017).


\bibitem[Se72]{Se72} {\sc J.-P. Serre}, {\it Propri\'{e}t\'{e}s galoisiennes des points d'ordre fini des courbes
   elliptiques}, Invent. Math., \textbf{15}(4) (1972), 259-331.
   




\bibitem[St75]{St75} {\sc R. Steinberg}, {\it On a theorem of Pittie}, Topology \textbf{14} (1975), 173-117.

\bibitem[Th87]{Th87} {\sc R.W. Thomason}, {\it Algebraic $K$-theory of group scheme actions}, Ann. of Math. Stud. \textbf{113} (1987), 539-563.

\bibitem[V94]{V94} {\sc M.-F. Vigneras}, {\it A propos d'une conjecture de Langlands modulaire}, in: {\it Finite reductive groups}, Ed. M. Cabanes, Prog. Math., vol. 141 (Birkhäuser, Basel, 1997).

\bibitem[V96]{V96} {\sc M.-F. Vigneras}, {\it Repr\'esentations $\ell$-modulaires d'un groupe r\'eductif p-adique}, Prog. Math., vol. 137 (Birkhäuser, Basel, 1996).

\bibitem[V04]{V04} {\sc M.-F. Vigneras} {\it Representations modulo $p$ of the $p$-adic group $GL(2,F)$}, Compositio Math. 140 (2004) 333-358.

\bibitem[V05]{V05} {\sc M.-F. Vigneras} {\it Pro-$p$-Iwahori Hecke ring and supersingular $\overline{\bbF}_p$-representations}, Math. Ann. \textbf{331} (2005), 523-556. + Erratum


\bibitem[V06]{V06} {\sc M.-F. Vigneras} {\it Algèbres de Hecke affines génériques}, Representation Theory \textbf{10} (2006), 1-20.

\bibitem[V14]{V14} {\sc M.-F. Vigneras}, {\it The pro-$p$-Iwahori Hecke algebra of a reductive $p$-adic group II}, Compositio Math. Muenster J. Math. \textbf{7} (2014), 363-379. + Erratum

    \bibitem[V15]{V15} {\sc M.-F. Vigneras}, {\it The pro-$p$-Iwahori Hecke algebra of a reductive $p$-adic group V (Parabolic induction)}, Pacific J. of Math. \textbf{279} (2015), Issue 1-2, 499-529.

\bibitem[V16]{V16} {\sc M.-F. Vigneras}, {\it The pro-$p$-Iwahori Hecke algebra of a reductive $p$-adic group I}, Compositio Math. \textbf{152} (2016), 693-753.

\bibitem[V17]{V17} {\sc M.-F. Vigneras}, {\it The pro-$p$-Iwahori Hecke algebra of a reductive $p$-adic group III (Spherical Hecke algebras and supersingular modules)}, Journal of the Institue of Mathematics of Jussieu \textbf{16} (2017), Issue 3, 571-608. + Erratum

\end{thebibliography}
\end{document}